\documentclass[twoside,twocolumn,english]{IEEEtran}
\usepackage[T1]{fontenc}
\usepackage{float}
\usepackage{amsthm}
\usepackage{amsbsy}
\usepackage{amstext}
\usepackage{graphicx}

\makeatletter

\providecommand{\tabularnewline}{\\}
\floatstyle{ruled}
\newfloat{algorithm}{tbp}{loa}
\providecommand{\algorithmname}{Algorithm}
\floatname{algorithm}{\protect\algorithmname}

\theoremstyle{plain}
\newtheorem{thm}{\protect\theoremname}
\theoremstyle{plain}
\newtheorem{lem}[thm]{\protect\lemmaname}
\theoremstyle{plain}
\newtheorem{prop}[thm]{\protect\propositionname}

\usepackage{amsmath}
\usepackage{amssymb}
\usepackage[numbers]{natbib}

\usepackage{graphicx}
\usepackage{subfigure}

\usepackage{pdfsync}
\usepackage{todonotes}
\makeatother

\usepackage{babel}
\providecommand{\lemmaname}{Lemma}
\providecommand{\propositionname}{Proposition}
\providecommand{\theoremname}{Theorem}

\begin{document}

\title{An Online Parallel and Distributed Algorithm for Recursive Estimation
of Sparse Signals}

\author{Yang Yang, Mengyi Zhang, Marius Pesavento, and Daniel P. Palomar%
\thanks{Y. Yang and M. Pesavento are with Communication Systems Group, Darmstadt
University of Technology, 64283 Darmstadt, Germany (Email: \{yang,pesavento\}@nt.tu-darmstadt.de).
Their work is supported by the Seventh Framework Programme for Research
of the European Commission under grant number: ADEL 619647. Yang's work was also supported by the Hong Kong RGC 16207814 research grant.%
}%
\thanks{M. Zhang is with Department of Computer Science and Engineering, The
Chinese University of Hong Kong, Hong Kong (Email: zhangmy@cse.cuhk.edu.hk). Her work was supported by the Hong Kong RGC 16207814 research grant.
}%
\thanks{D. P. Palomar is with Department of Electronic and Computer Engineering,
The Hong Kong University of Science and Technology, Hong Kong (Email:
palomar@ust.hk). His work is supported by the Hong Kong RGC 16207814 research grant.
}
\thanks{Part of this work has been presented at The Asilomar Conference on Signals, Systems, and Computers, Nov. 2014.
}}
\maketitle
\begin{abstract}
In this paper, we consider a recursive estimation problem for linear
regression where the signal to be estimated admits a sparse representation
and measurement samples are only sequentially available. We propose
a convergent parallel estimation scheme that consists in solving a
sequence of $\ell_{1}$-regularized least-square problems approximately.
The proposed scheme is novel in three aspects: i) all elements of
the unknown vector variable are updated in parallel at each time instance,
and convergence speed is much faster than state-of-the-art schemes
which update the elements sequentially; ii) both the update direction
and stepsize of each element have simple closed-form expressions,
so the algorithm is suitable for online (real-time) implementation;
and iii) the stepsize is designed to accelerate the convergence but
it does not suffer from the common trouble of parameter tuning in
literature. Both centralized and distributed implementation schemes
are discussed. The attractive features of the proposed algorithm are
also numerically consolidated.
\end{abstract}

\section{Introduction}

Signal estimation has been a fundamental problem in a number of scenarios,
such as wireless sensor networks (WSN) and cognitive radio (CR). WSN
has received a lot of attention and is found to be useful in diverse
disciplines such as environmental monitoring, smart grid, and wireless
communications \cite{Barbarossa2013}. CR appears as an enabling technique
for flexible and efficient use of the radio spectrum \cite{Haykin2005a,Mitola1999},
since it allows the unlicensed secondary users (SUs) to access the
spectrum provided that the licensed primary users (PUs) are idle,
and/or the interference generated by the SUs to the PUs is below a
certain level that is tolerable for the PUs \cite{Zhang2010a,Yang2013JSAC}.

One in CR systems is the ability to obtain a precise estimate of the
PUs' power distribution map so that the SUs can avoid the areas in
which the PUs are actively transmitting. This is usually realized
through the estimation of the position, transmit status, and/or transmit
power of PUs \cite{Haykin2009,Kim2011a,Zeng2011,Mehanna2013}, and
such an estimation is typically obtained based on the minimum mean-square-error
(MMSE) criterion \cite{Estimation_Kay,Sayed2008Adaptive,Mateos2009,Mateos2010,Angelosante2010,Zeng2011,Barbarossa2013}.

The MMSE approach involves the calculation of the expectation of a
squared $\ell_{2}$-norm function that depends on the so-called regression
vector and measurement output, both of which are random variables.
This is essentially a stochastic optimization problem, but when the
statistics of these random variables are unknown, it is impossible
to calculate the expectation analytically. An alternative is to use
the sample average function, constructed from the sequentially available
measurements, as an approximation of the expectation, and this leads
to the well-known recursive least-square (RLS) algorithm \cite{Sayed2008Adaptive,Mateos2009,Mateos2010,Barbarossa2013}.
As the measurements are available sequentially, at each time instance
of the RLS algorithm, an LS problem has to be solved, which furthermore
admits a closed-form solution and thus can efficiently be computed.
More details can be found in standard textbooks such as \cite{Estimation_Kay,Sayed2008Adaptive}.

In practice, the signal to be estimated may be sparse in nature \cite{Angelosante2010,Kim2011a,Kopsinis2011,Zeng2011,Barbarossa2013}.
In a recent attempt to apply the RLS approach to estimate a sparse
signal, a regularization function in terms of $\ell_{1}$-norm was
incorporated into the LS function to encourage sparse estimates \cite{Angelosante2010,Barbarossa2013},
leading to an $\ell_{1}$-regularized LS problem which has the form
of the least-absolute shrinkage and selection operator (LASSO) \cite{Tibshirani1996}.
Then in the recursive estimation of a sparse signal, the only difference
from standard RLS is that at each time instance, instead of solving
an LS problem as in RLS, an $\ell_{1}$-regularized LS problem in
the form of LASSO is solved \cite{Barbarossa2013}.

However, a closed-form solution to the $\ell_{1}$-regularized LS
problem no longer exists because of the $\ell_{1}$-norm regularization
function and the problem can only be solved iteratively. As a matter
of fact, iterative algorithms to solve the $\ell_{1}$-regularized
LS problems have been the center of extensive research in recent years
and a number of solvers have been developed, e.g., $\texttt{GP}$
\cite{Figueiredo2007}, $\texttt{l1\_ls}$ \cite{Kim2007}, $\texttt{FISTA}$
\cite{Beck2009}, $\texttt{ADMM}$ \cite{Goldstein2009}, and $\texttt{FLEXA}$
\cite{Scutari_BigData}. Since the measurements are sequentially available,
and with each measurement, a new $\ell_{1}$-regularized LS problem
is formed and solved, the overall complexity of using solvers for
the whole sequence of $\ell_{1}$-regularized LS problems is no longer
affordable. If the environment is furthermore fast changing, this
method is not even real-time applicable because new samples may have
already arrived before the old $\ell_{1}$-regularized LS problem
is solved.

To make the estimation scheme suitable for online (real-time) implementation,
a\emph{ }sequential algorithm was proposed in \cite{Angelosante2010},
in which the $\ell_{1}$-regularized LS problem at each time instance
is solved only approximately. In particular, at each time instance,
the $\ell_{1}$-regularized LS problem is solved with respect to (w.r.t.)
only a single element of the unknown vector variable (instead of \emph{all
}elements as in a solver) while remaining elements are fixed, and
the element is updated in closed-form based on the so-called soft-thresholding
operator \cite{Beck2009}. After a new sample arrives, a new $\ell_{1}$-regularized
LS problem is formed and solved w.r.t. the next element while remaining
elements are fixed. This sequential update rule is known in literature
as block coordinate descent method \cite{Bertsekas}. To our best
knowledge, \cite{Angelosante2010} is the only work on online algorithms
for recursive estimation of sparse signals.

Intuitively, since only a single element is updated at each time instance,
the online algorithm proposed in \cite{Angelosante2010} sometimes
suffers from slow convergence, especially when the signal has a large
dimension while large dimension of sparse signals is universal in
practice. It is tempting to use the parallel algorithm proposed in
\cite{Scutarib,Scutari_BigData}, but it works for deterministic optimization
problems only and may not converge for the stochastic optimization
problem at hand. Besides, its convergence speed heavily depends on
the stepsize. Typical stepsizes are Armijo-like successive line search,
constant stepsize, and diminishing stepsize. The former two suffer
from high complexity and slow convergence \cite[Remark 4]{Scutari_BigData},
while the decay rate of the diminishing stepsize is very difficult
to choose: on the one hand, a slowly decaying stepsize is preferable
to make notable progress and to achieve satisfactory convergence speed;
on the other hand, theoretical convergence is guaranteed only when
the stepsizes decays fast enough. It is a difficult task on its own
to find the decay rate that gives a good trade-off.

A recent work on parallel algorithms for stochastic optimization is
\cite{Yang_stochastic}. However, the algorithms proposed in \cite{Yang_stochastic}
are not applicable for the recursive estimation of sparse signals.
This is because the regularization function in \cite{Yang_stochastic}
must be strongly convex and differentiable while the regularization
gain must be lower bounded by some positive constant so that convergence
can be achieved, but the regularization function in terms of $\ell_{1}$-norm
in this paper is convex (but not strongly convex) and nonsmooth while
the regularization gain is decreasing to 0.

In this paper, we propose an online parallel algorithm with provable\emph{
}convergence for recursive estimation of sparse signals. In particular,
our contributions are as follows:

1) At each time instance, the $\ell_{1}$-regularized LS problem is
solved approximately and all elements are updated in parallel, so
the convergence speed is greatly enhanced compared with \cite{Angelosante2010}.
As a nontrivial extension of \cite{Angelosante2010} from sequential
update to parallel update and \cite{Scutarib,Scutari_BigData} from
deterministic optimization problems to stochastic optimization problems,
the convergence of the proposed algorithm is established.

2) The proposed stepsize is based on the so-called minimization rule
(also known as exact line search) and its benefits are twofold: firstly,
it guarantees the convergence of the proposed algorithm, which may
however diverge under other stepsize rules; secondly, notable progress
is achieved after each variable update and the trouble of parameter
tuning in \cite{Scutarib,Scutari_BigData} is saved. Besides, both
the update direction and stepsize of each element have a simple closed-form
expression, so the algorithm is fast to converge and suitable for
online implementation.

3) When implemented in a distributed manner, for example in CR networks,
the proposed algorithm has a much smaller signaling overhead than
in state-of-the-art techniques \cite{Barbarossa2013,Kim2011a}. Besides,
the estimates of different SUs are always the same and they are based
on the global information of all SUs. Compared with consensus-based
distributed implementations \cite{Zeng2011,Kim2011a,Barbarossa2013}
where each SU makes individual estimate decision mainly based on his
own local information and all SUs converge to the same estimate only
asymptotically, the proposed approach can better protect the Quality-of-Service
(QoS) of PUs because it eliminates the possibility that the estimates
maintained by different SUs lead to conflicting interests of the PUs
(i.e., some may correctly detect the presence of PUs but some may
not in consensus-based algorithms).

The rest of the paper is organized as follows. In Section \ref{sec:System-Model}
we introduce the system model and formulate the recursive estimation
problem. The online parallel algorithm is proposed in Section \ref{sec:algorithm},
and its implementations and extensions are discussed in Section \ref{sec:Implementation-and-Extensions}.
The performance of the proposed algorithm is evaluated numerically
in Section \ref{sec:Numerical-Results} and finally concluding remarks
are drawn in Section \ref{sec:Concluding-Remarks}.

\emph{Notation: }We use $x$, $\mathbf{x}$ and $\mathbf{X}$ to denote
scalar, vector and matrix, respectively. $X_{jk}$ is the $(j,k)$-th
element of $\mathbf{X}$; $x_{k}$ and $x_{j,k}$ is the $k$-th element
of $\mathbf{x}$ and $\mathbf{x}_{j}$, respectively, and $\mathbf{x}=(x_{k})_{k=1}^{K}$
and $\mathbf{x}_{j}=(x_{j,k})_{k=1}^{K}$. $\mathbf{d(X)}$ is a vector
that consists of the diagonal elements of $\mathbf{X}$. $\mathbf{x}\circ\mathbf{y}$
denotes the Hadamard product between $\mathbf{x}$ and $\mathbf{y}$.
$[\mathbf{x}]_{\mathbf{a}}^{\mathbf{b}}$ denotes the element-wise
projection of $\mathbf{x}$ onto $[\mathbf{a,b}]$: $[\mathbf{x}]_{\mathbf{a}}^{\mathbf{b}}\triangleq\max(\min(\mathbf{x},\mathbf{b}),\mathbf{a})$,
and $\left[\mathbf{x}\right]^{+}$ denotes the element-wise projection
of $\mathbf{x}$ onto the nonnegative orthant: $\left[\mathbf{x}\right]^{+}\triangleq\max(\mathbf{x},\mathbf{0})$.
$\mathbf{X}^{\dagger}$ denotes the Moore-Penrose inverse of $\mathbf{X}$.

\section{\label{sec:System-Model}System Model and Problem Formulation}

Suppose $\mathbf{x}^{\star}=(x_{k}^{\star})_{k=1}^{K}\in\mathbb{R}^{K}$
is a deterministic sparse signal to be estimated based on the the
measurement $y_{n}\in\mathbb{R}$, and they are connected through
a linear regression model:
\begin{equation}
y_{n}=\mathbf{g}_{n}^{T}\mathbf{x}^{\star}+v_{n},\quad n=1,\ldots,N,\label{eq:linear-model}
\end{equation}
where $N$ is the number of measurements at any time instance. The
regression vector $\mathbf{g}_{n}=(g_{n,k})_{k=1}^{K}\in\mathbb{R}^{K}$
is assumed to be known, and $v_{n}\in\mathbb{R}$ is the additive
estimation noise. Throughout the paper, we make the following assumptions
on $\mathbf{g}_{n}$ and $v_{n}$ for $n=1,\ldots,N$:

\begin{enumerate}

\item[(A1)] $\mathbf{g}_{n}$ are independently and identically distributed
(i.i.d.) random variables with a bounded positive definite covariance
matrix;

\item[(A2)] $v_{n}$ are i.i.d. random variables with zero mean and
bounded variance, and are uncorrelated with $\mathbf{g}_{n}$.

\end{enumerate}

Given the linear model in (\ref{eq:linear-model}), the problem is
to estimate $\mathbf{x}^{\star}$ from the set of regression vectors
and measurements $\left\{ \mathbf{g}_{n},y_{n}\right\} _{n=1}^{N}$.
Since both the regression vector $\mathbf{g}_{n}$ and estimation
noise $v_{n}$ are random variables, the measurement $y_{n}$ is also
random. A fundamental approach to estimate $\mathbf{x}^{\star}$ is
based on the MMSE criterion, which has a solid root in adaptive filter
theory \cite{Sayed2008Adaptive,Estimation_Kay}. To improve the estimation
precision, all available measurements $\left\{ \mathbf{g}_{n},y_{n}\right\} _{n=1}^{N}$
are exploited to form a cooperative estimation problem which consists
in finding the variable that minimizes the mean-square-error \cite{Quan2008,Barbarossa2013,Zeng2011}:
\begin{align}
\mathbf{x}^{\star} & =\underset{\mathbf{x}=(x_{k})_{k=1}^{K}}{\arg\min}\;\mathbb{E}\left[\sum_{n=1}^{N}\left(y_{n}-\mathbf{g}_{n}^{T}\mathbf{x}\right)^{2}\right]\label{eq:mmse}\\
 & =\underset{\mathbf{x}}{\arg\min}\;\frac{1}{2}\mathbf{x}^{T}\mathbf{G}\mathbf{x}-\mathbf{b}^{T}\mathbf{x},\nonumber
\end{align}
where $\mathbf{G}\triangleq\sum_{n=1}^{N}\mathbb{E}\left[\mathbf{g}_{n}\mathbf{g}_{n}^{T}\right]$
and $\mathbf{b}\triangleq\sum_{n=1}^{N}\mathbb{E}\left[y_{n}\mathbf{g}_{n}\right]$,
and the expectation is taken over $\{\mathbf{g}_{n},y_{n}\}_{n=1}^{N}$.

In practice, the statistics of $\{\mathbf{g}_{n},y_{n}\}_{n=1}^{N}$
are often not available to compute $\mathbf{G}$ and $\mathbf{b}$
analytically. In fact, the absence of statistical information is a
general rule rather than an exception. It is a common approach to
approximate the expectation in (\ref{eq:mmse}) by the sample average
constructed from the samples $\{\mathbf{g}_{n}^{(\tau)},y_{n}^{(\tau)}\}_{\tau=1}^{t}$
sequentially available up to time $t$ \cite{Sayed2008Adaptive}:\begin{subequations}\label{eq:rls}
\begin{align}
\mathbf{x}_{\textrm{rls}}^{(t)} & \triangleq\underset{\mathbf{x}}{\arg\min}\;\frac{1}{2}\mathbf{x}^{T}\mathbf{G}^{(t)}\mathbf{x}-(\mathbf{b}^{(t)})^{T}\mathbf{x}\label{eq:rls-1}\\
 & =\mathbf{G}^{(t)\dagger}\mathbf{b}^{(t)},\label{eq:rls-2}
\end{align}
\end{subequations}where $\mathbf{G}^{(t)}$ and $\mathbf{b}^{(t)}$
is the sample average of $\mathbf{G}$ and $\mathbf{b}$, respectively:
\begin{equation}
\mathbf{G}^{(t)}\triangleq\frac{1}{t}\sum_{\tau=1}^{t}\sum_{n=1}^{N}\mathbf{g}_{n}^{(\tau)}(\mathbf{g}_{n}^{(\tau)})^{T},\;\mathbf{b}^{(t)}\triangleq\frac{1}{t}\sum_{\tau=1}^{t}\sum_{n=1}^{N}y_{n}^{(\tau)}\mathbf{g}_{n}^{(\tau)},\label{eq:A-and-b}
\end{equation}
and $\mathbf{A}^{\dagger}$ is the Moore-Penrose pseudo-inverse of
$\mathbf{A}$. In literature, (\ref{eq:rls}) is known as recursive
least square (RLS), as indicated by the subscript ``rls'', and $\mathbf{x}_{\textrm{rls}}^{(t)}$
can be computed efficiently in closed-form, cf. (\ref{eq:rls-2}).

In many practical applications, the unknown signal $\mathbf{x}^{\star}$
is sparse by nature or by design, but $\mathbf{x}_{\textrm{rls}}^{(t)}$
given by (\ref{eq:rls}) is not necessarily sparse when $t$ is finite
\cite{Tibshirani1996,Kim2007}. To overcome this shortcoming, a sparsity
encouraging function in terms of $\ell_{1}$-norm is incorporated
into the sample average function in (\ref{eq:rls}), leading to the
following $\ell_{1}$-regularized sample average function at any time
instance $t=1,2,\ldots$ \cite{Angelosante2010,Kim2011a,Barbarossa2013}:
\begin{equation}
L^{(t)}(\mathbf{x})\triangleq\frac{1}{2}\mathbf{x}^{T}\mathbf{G}^{(t)}\mathbf{x}-(\mathbf{b}^{(t)})^{T}\mathbf{x}+\mu^{(t)}\left\Vert \mathbf{x}\right\Vert _{1},\label{eq:L^t(x)}
\end{equation}
where $\mu^{(t)}>0$. Define $\mathbf{x}_{\textrm{lasso}}^{(t)}$
as the minimizing variable of $L^{(t)}(\mathbf{x})$:
\begin{equation}
\mathbf{x}_{\textrm{lasso}}^{(t)}=\underset{\mathbf{x}}{\arg\min}\; L^{(t)}(\mathbf{x}),\quad t=1,2,\ldots,\label{eq:l1-RLS}
\end{equation}
In literature, problem (\ref{eq:l1-RLS}) for any fixed $t$ is known
as the \emph{least-absolute shrinkage and selection operator }(LASSO)
\cite{Tibshirani1996,Kim2007} (as indicated by the subscript ``lasso''
in (\ref{eq:l1-RLS})). Note that in batch processing \cite{Tibshirani1996,Kim2007},
problem (\ref{eq:l1-RLS}) is solved only once when a certain number
of measurements are collected (so $t$ is equal to the number of measurements),
while in the recursive estimation of $\mathbf{x}^{\star}$, the measurements
are sequentially available (so $t$ is increasing) and (\ref{eq:l1-RLS})
is solved repeatedly at each time instance $t=1,2,\ldots$

The advantage of (\ref{eq:l1-RLS}) over (\ref{eq:mmse}), whose objective
function is stochastic and whose calculation depends on unknown parameters
$\mathbf{G}$ and $\mathbf{b}$, is that (\ref{eq:l1-RLS}) is a sequence
of deterministic optimization problems whose theoretical and algorithmic
properties have been extensively investigated and widely understood.
A natural question arises in this context: is (\ref{eq:l1-RLS}) equivalent
to (\ref{eq:mmse}) in the sense that $\mathbf{x}_{\textrm{lasso}}^{(t)}$
is a strongly consistent estimator of $\mathbf{x}^{\star}$, i.e.,
$\lim_{t\rightarrow\infty}\mathbf{x}_{\textrm{lasso}}^{(t)}=\mathbf{x}^{\star}$
with probability 1? The connection between $\mathbf{x}_{\textrm{lasso}}^{(t)}$
in (\ref{eq:l1-RLS}) and the unknown variable $\mathbf{x}^{\star}$
is given in the following lemma \cite{Angelosante2010}.
\begin{lem}
\label{lem:connection}Suppose Assumptions (A1)-(A2) as well as the
following assumption are satisfied for (\ref{eq:l1-RLS}):

\begin{enumerate}

\item[\em{(A3)}] $\left\{ \mu^{(t)}\right\} $ is a positive sequence
converging to $0$, i.e., $\mu^{(t)}>0$ and $\lim_{t\rightarrow\infty}\mu^{(t)}=0$.

\end{enumerate}Then $\lim_{t\rightarrow\infty}\mathbf{x}_{\mathrm{lasso}}^{(t)}=\mathbf{x}^{\star}$
with probability 1.
\end{lem}
An example of $\mu^{(t)}$ satisfying Assumption (A3) is $\mu^{(t)}=\alpha/t^{\beta}$
with $\alpha>0$ and $\beta>0$. Typical choices of $\beta$ are $\beta=1$
and $\beta=0.5$ \cite{Angelosante2010}.

Lemma \ref{lem:connection} not only states the connection between
$\mathbf{x}_{\textrm{lasso}}^{(t)}$ and $\mathbf{x}^{\star}$ from
a theoretical perspective, but also suggests a simple algorithmic
solution for problem (\ref{eq:mmse}): $\mathbf{x}^{\star}$ can be
estimated by solving a sequence of deterministic optimization problems
(\ref{eq:l1-RLS}), one for each time instance $t=1,2,\ldots$. However,
different from RLS in which each update has a closed-form expression,
cf. (\ref{eq:rls-2}), problem (\ref{eq:l1-RLS}) does not have a
closed-form solution and it can only be solved numerically by iterative
algorithm such as $\texttt{GP}$ \cite{Figueiredo2007}, $\texttt{l1\_ls}$
\cite{Kim2007}, $\texttt{FISTA}$ \cite{Beck2009}, $\texttt{ADMM}$
\cite{Goldstein2009}, and $\texttt{FLEXA}$ \cite{Scutari_BigData}.
As a result, solving (\ref{eq:l1-RLS}) repeatedly at each time instance
$t=1,2,\ldots$ is neither computationally practical nor real-time
applicable. The aim of the following sections is to develop an algorithm
that enjoys easy implementation and fast convergence.

\section{\label{sec:algorithm}The Online Parallel Algorithm}

The LASSO problem in (\ref{eq:l1-RLS}) is convex, but the objective
function is nondifferentiable and it cannot be minimized in closed-form,
so solving (\ref{eq:l1-RLS}) completely w.r.t. all elements of $\mathbf{x}$
by a solver at each time instance $t=1,2,\ldots$ is neither computationally
practical nor suitable for online implementation. To reduce the complexity
of the variable update, an algorithm based on inexact optimization
is proposed in \cite{Angelosante2010}: at time instance $t$, only
a single element $x_{k}$ with $k=\textrm{mod}(t-1,K)+1$ is updated
by its so-called best response, i.e., $L^{(t)}(\mathbf{x})$ is minimized
w.r.t. $x_{k}$ only: $x_{k}^{(t+1)}=\arg\min L^{(t)}(x_{k},\mathbf{x}_{-k}^{(t)})$
with $\mathbf{x}_{-k}\triangleq(x_{j})_{j\neq k}$, which can be solved
in closed-form, while the remaining elements $\left\{ x_{j}\right\} _{j\neq k}$
remain unchanged, i.e., $\mathbf{x}_{-k}^{(t+1)}=\mathbf{x}_{-k}^{(t)}$.
At the next instance $t+1$, a new sample average function $L^{(t+1)}(\mathbf{x})$
is formed with newly arriving samples, and the $(k+1)$-th element,
$x_{k+1}$, is updated by minimizing $L^{(t+1)}(\mathbf{x})$ w.r.t.
$x_{k+1}$ only, while the remaining elements again are fixed. Although
easy to implement, sequential updating schemes update only a single
element at each time instance and they sometimes suffer from slow
convergence when the number of elements $K$ is large.

To overcome the slow convergence of the sequential update, we propose
an online parallel update scheme, with provable convergence, in which
(\ref{eq:l1-RLS}) is solved \emph{approximately }by simultaneously
updating all elements only once based on their individual best response.
Given the current estimate $\mathbf{x}^{(t)}$ which is available
before the $t$-th sample arrives%
\footnote{$\mathbf{x}^{(1)}$ could be arbitrarily chosen, e.g., $\mathbf{x}^{(1)}=\mathbf{0}$.%
}, the next estimate $\mathbf{x}^{(t+1)}$ is determined based on all
the samples collected up to instance $t$ in a three-step procedure
as described next.

 \newcounter{MYtempeqncnt} \begin{figure*}[t] \normalsize \setcounter{MYtempeqncnt}{\value{equation}} \setcounter{equation}{15} \vspace*{4pt} \begin{equation} \label{eq:stepsize-2}  \gamma^{(t)}=\left[-\frac{(\mathbf{G}^{(t)}\mathbf{x}^{(t)}-\mathbf{b}^{(t)})^T (\hat{\mathbf{x}}^{(t)}-\mathbf{x}^{(t)})+\mu^{(t)}(\bigl\Vert\hat{\mathbf{x}}^{(t)}\bigr\Vert_{1}-\bigl\Vert\mathbf{x}^{(t)}\bigr\Vert_{1})}{(\hat{\mathbf{x}}^{(t)}-\mathbf{x}^{(t)})^T\mathbf{G}^{(t)}(\hat{\mathbf{x}}^{(t)}-\mathbf{x}^{(t)})}\right]^{1}_{0} \end{equation} \setcounter{equation}{\value{MYtempeqncnt}} \hrulefill  \end{figure*}

\textbf{Step 1 (Update Direction): }In this step, all elements of
$\mathbf{x}$ are updated \emph{in parallel} and the update direction
of $\mathbf{x}$ at $\mathbf{x}=\mathbf{x}^{(t)}$, denoted as $\hat{\mathbf{x}}^{(t)}-\mathbf{x}^{(t)}$,
is determined based on the best-response $\hat{\mathbf{x}}^{(t)}$.
For each element of $\mathbf{x}$, say $x_{k}$, its best response
at $\mathbf{x}=\mathbf{x}^{(t)}$ is given by:
\begin{equation}
\hat{x}_{k}^{(t)}\triangleq\underset{x_{k}}{\arg\min}\;\Bigl\{ L^{(t)}(x_{k},\mathbf{x}_{-k}^{(t)})+\frac{1}{2}c_{k}^{(t)}(x_{k}-x_{k}^{(t)})^{2}\Bigr\},\quad\forall\, k,\label{eq:x-hat-definition}
\end{equation}
where $\mathbf{x}_{-k}\triangleq\{x_{j}\}_{j\neq k}$ and it is fixed
to their values of the preceding time instance $\mathbf{x}_{-k}=\mathbf{x}_{-k}^{(t)}$.
An additional quadratic proximal term with $c_{k}^{(t)}>0$ is included
in (\ref{eq:x-hat-definition}) for numerical simplicity and stability
\cite{Bertsekas,Scutarib}, because it plays an important role in
the convergence analysis of the proposed algorithm; conceptually it
is a penalty (with variable weight $c_{k}^{(t)}$) for moving away
from the current estimate $x_{k}^{(t)}$.

After substituting (\ref{eq:L^t(x)}) into (\ref{eq:x-hat-definition}),
the best-response in (\ref{eq:x-hat-definition}) can be expressed
in closed-form:
\begin{align}
\hat{x}_{k}^{(t)} & =\underset{x_{k}}{\arg\min}\left\{ \negthickspace\negthickspace\begin{array}{l}
\frac{1}{2}G_{kk}^{(t)}x_{k}^{2}-r_{k}^{(t)}\cdot x_{k}\smallskip\\
\;+\mu^{(t)}|x_{k}|+\frac{1}{2}c_{k}^{(t)}(x_{k}-x_{k}^{(t)})^{2}
\end{array}\negthickspace\negthickspace\right\} \nonumber \\
 & =\frac{\mathcal{S}_{\mu^{(t)}}(r_{k}^{(t)}+c_{k}^{(t)}x_{k}^{(t)})}{G_{kk}^{(t)}+c_{k}^{(t)}},\quad k=1,\ldots,K,\label{eq:x-hat-scalar}
\end{align}
where
\begin{equation}
r_{k}^{(t)}\triangleq b_{k}^{(t)}-\sum_{j\neq k}G_{kj}^{(t)}x_{j}^{(t)},\label{eq:x-hat-r_t}
\end{equation}
and
\[
\mathcal{S}_{a}(b)\triangleq(b-a)^{+}-(-b-a)^{+}
\]
is the well-known soft-thresholding operator \cite{Beck2009,Boyd2010}.
From the definition of $\mathbf{G}^{(t)}$ in (\ref{eq:A-and-b}),
$\mathbf{G}^{(t)}\succeq\mathbf{0}$ and $G_{kk}^{(t)}\geq0$ for
all $k$, so the division in (\ref{eq:x-hat-scalar}) is well-defined.

Given the update direction $\hat{\mathbf{x}}^{(t)}-\mathbf{x}^{(t)}$,
an intermediate update vector $\tilde{\mathbf{x}}^{(t)}(\gamma)$
is defined:
\begin{equation}
\tilde{\mathbf{x}}^{(t)}(\gamma)=\mathbf{x}^{(t)}+\gamma(\hat{\mathbf{x}}^{(t)}-\mathbf{x}^{(t)}),\label{eq:x_tilde}
\end{equation}
where $\hat{\mathbf{x}}^{(t)}=(\hat{x}_{k}^{(t)})_{k=1}^{K}$ and
$\gamma\in[0,1]$ is the stepsize. The update direction $\hat{\mathbf{x}}^{(t)}-\mathbf{x}^{(t)}$
is a descent direction of $L^{(t)}(\mathbf{x})$ in the sense specified
by the following proposition.
\begin{prop}
[\bf{Descent Direction}]\label{prop:descent-direction}For $\hat{\mathbf{x}}^{(t)}=(\hat{x}_{k}^{(t)})_{k=1}^{K}$
given in (\ref{eq:x-hat-scalar}) and the update direction $\hat{\mathbf{x}}^{(t)}-\mathbf{x}^{(t)}$,
the following holds for any $\gamma\in[0,1]$:
\begin{align}
L^{(t)}(\tilde{\mathbf{x}}^{(t)}(\gamma))-L^{(t)}(\mathbf{x}^{(t)})\nonumber \\
\leq-\gamma\Bigl(c_{\min}^{(t)}-{\textstyle \frac{1}{2}}\lambda_{\max} & (\mathbf{G}^{(t)})\gamma\Bigr)\bigl\Vert\hat{\mathbf{x}}^{(t)}-\mathbf{x}^{(t)}\bigr\Vert_{2}^{2},\label{eq:descent}
\end{align}
where $c_{\min}^{(t)}\triangleq\min_{k}\left\{ G_{kk}^{(t)}+c_{k}^{(t)}\right\} >0$.\end{prop}
\begin{IEEEproof}
The proof follows the general line of arguments in \cite[Prop. 8(c)]{Scutari_BigData}
and is thus omitted here.
\end{IEEEproof}
\textbf{Step 2 (Stepsize): }In this step, the stepsize $\gamma$ in
(\ref{eq:x_tilde}) is determined so that fast convergence is observed.
It is easy to see from (\ref{eq:descent}) that for sufficiently small
$\gamma$, the right hand side of (\ref{eq:descent}) becomes negative
and $L^{(t)}(\mathbf{x})$ decreases as compared to $\mathbf{x}=\mathbf{x}^{(t)}$.
Thus, to minimize $L^{(t)}(\mathbf{x})$, a natural choice of the
stepsize rule is the so-called ``minimization rule'' \cite[Sec. 2.2.1]{bertsekas1999nonlinear}
(also known as the ``exact line search'' \cite[Sec. 9.2]{boyd2004convex}),
which is the stepsize, denoted as $\gamma_{\textrm{opt}}^{(t)}$,
that decreases $L^{(t)}(\mathbf{x})$ to the largest extent along
the direction $\hat{\mathbf{x}}^{(t)}-\mathbf{x}^{(t)}$ at $\mathbf{x}=\mathbf{x}^{(t)}$:
\begin{align}
\gamma_{\textrm{opt}}^{(t)} & =\underset{0\leq\gamma\leq1}{\arg\min}\;\bigl\{ L^{(t)}(\tilde{\mathbf{x}}^{(t)}(\gamma))-L^{(t)}(\mathbf{x}^{(t)})\bigr\}\nonumber \\
 & =\underset{0\leq\gamma\leq1}{\arg\min}\;\bigl\{ L^{(t)}(\mathbf{x}^{(t)}+\gamma(\hat{\mathbf{x}}^{(t)}-\mathbf{x}^{(t)}))-L^{(t)}(\mathbf{x}^{(t)})\bigr\}\nonumber \\
 & =\underset{0\leq\gamma\leq1}{\arg\min}\negthickspace\left\{ \negthickspace\negthickspace\negthickspace\begin{array}{l}
\frac{1}{2}(\hat{\mathbf{x}}^{(t)}-\mathbf{x}^{(t)})^{T}\mathbf{G}^{(t)}(\hat{\mathbf{x}}^{(t)}-\mathbf{x}^{(t)})\cdot\gamma^{2}\smallskip\\
\;+(\mathbf{G}^{(t)}\mathbf{x}^{(t)}-\mathbf{b}^{(t)})^{T}(\hat{\mathbf{x}}^{(t)}-\mathbf{x}^{(t)})\cdot\gamma\smallskip\\
\;+\mu^{(t)}(\left\Vert \mathbf{x}^{(t)}+\gamma(\hat{\mathbf{x}}^{(t)}-\mathbf{x}^{(t)})\right\Vert _{1}-\left\Vert \mathbf{x}^{(t)}\right\Vert _{1})
\end{array}\negthickspace\negthickspace\negthickspace\right\} \negthickspace.\label{eq:minimization-rule-1}
\end{align}
Therefore by definition of $\gamma_{\textrm{opt}}^{(t)}$ we have
for any $\gamma\in[0,1]$:
\begin{equation}
L^{(t)}(\mathbf{x}^{(t)}+\gamma_{\textrm{opt}}^{(t)}(\hat{\mathbf{x}}^{(t)}-\mathbf{x}^{(t)}))\leq L^{(t)}(\mathbf{x}^{(t)}+\gamma(\hat{\mathbf{x}}^{(t)}-\mathbf{x}^{(t)})).\label{eq:minimization-rule-2}
\end{equation}
The difficulty with the standard minimization rule (\ref{eq:minimization-rule-1})
is the complexity of solving the optimization problem in (\ref{eq:minimization-rule-1}),
since the presence of the $\ell_{1}$-norm makes it impossible to
find a closed-form solution and the problem in (\ref{eq:minimization-rule-1})
can only be solved numerically by a solver such as $\texttt{SeDuMi}$
\cite{Sturm1999}.

To obtain a stepsize with a good trade off between convergence speed
and computational complexity, we propose a \emph{simplified} minimization
rule which yields fast convergence but can be computed at a low complexity.
Firstly it follows from the convexity of norm functions that for any
$\gamma\in[0,1]$:\begin{subequations}\label{eq:triangle-inequality}
\begin{align}
 & \mu^{(t)}(\bigl\Vert\mathbf{x}^{(t)}+\gamma(\hat{\mathbf{x}}^{(t)}-\mathbf{x}^{(t)})\bigr\Vert_{1}-\bigl\Vert\mathbf{x}^{(t)}\bigr\Vert_{1})\nonumber \\
=\; & \mu^{(t)}\bigl\Vert(1-\gamma)\mathbf{x}^{(t)}+\gamma\hat{\mathbf{x}}^{(t)}\bigr\Vert_{1}-\mu^{(t)}\bigl\Vert\mathbf{x}^{(t)}\bigr\Vert_{1}\nonumber \\
\leq\; & (1-\gamma)\mu^{(t)}\bigl\Vert\mathbf{x}^{(t)}\bigr\Vert_{1}+\gamma\mu^{(t)}\bigl\Vert\hat{\mathbf{x}}^{(t)}\bigr\Vert_{1}-\mu^{(t)}\bigl\Vert\mathbf{x}^{(t)}\bigr\Vert_{1}\label{eq:triangle-inequality-1}\\
=\; & \mu^{(t)}(\bigl\Vert\hat{\mathbf{x}}^{(t)}\bigr\Vert_{1}-\bigl\Vert\mathbf{x}^{(t)}\bigr\Vert_{1})\cdot\gamma.\label{eq:triangle-inequality-2}
\end{align}
\end{subequations}The right hand side of (\ref{eq:triangle-inequality-2})
is linear in $\gamma$, and equality is achieved in (\ref{eq:triangle-inequality-1})
either when $\gamma=0$ or $\gamma=1$.

In the proposed simplified minimization rule, instead of directly
minimizing $L^{(t)}(\tilde{\mathbf{x}}^{(t)}(\gamma))-L^{(t)}(\mathbf{x}^{(t)})$
over $\gamma$, its upper bound based on (\ref{eq:triangle-inequality})
is minimized and $\gamma^{(t)}$ is given by
\begin{equation}
\gamma^{(t)}\triangleq\underset{0\leq\gamma\leq1}{\arg\min}\left\{ \negthickspace\negthickspace\begin{array}{l}
\frac{1}{2}(\hat{\mathbf{x}}^{(t)}-\mathbf{x}^{(t)})^{T}\mathbf{G}^{(t)}(\hat{\mathbf{x}}^{(t)}-\mathbf{x}^{(t)})\cdot\gamma^{2}\smallskip\\
\;+(\mathbf{G}^{(t)}\mathbf{x}^{(t)}-\mathbf{b}^{(t)})^{T}(\hat{\mathbf{x}}^{(t)}-\mathbf{x}^{(t)})\cdot\gamma\smallskip\\
\;+\mu^{(t)}(\left\Vert \hat{\mathbf{x}}^{(t)}\right\Vert _{1}-\left\Vert \mathbf{x}^{(t)}\right\Vert _{1})\cdot\gamma\smallskip
\end{array}\negthickspace\negthickspace\right\} ,\label{eq:stepsize-1}
\end{equation}
The scalar problem in (\ref{eq:stepsize-1}) consists in a convex
quadratic objective function along with a bound constraint and it
has a closed-form solution, given by (\ref{eq:stepsize-2}) at the
top of the next page, where $[x]_{0}^{1}\triangleq\min(\max(x,0),1)$
denotes the projection of $x$ onto $[0,1]$, and obtained by projecting
the unconstrained optimal variable of the convex quadratic scalar
problem in (\ref{eq:stepsize-1}) onto the interval $[0,1]$.

The advantage of minimizing the upper bound function of $L^{(t)}(\tilde{\mathbf{x}}^{(t)}(\gamma))$
in (\ref{eq:stepsize-1}) is that the optimal $\gamma$, denoted as
$\gamma^{(t)}$, always has a closed-form expression, cf. (\ref{eq:stepsize-2}).
At the same time, it also yields a decrease in $L^{(t)}(\mathbf{x})$
at $\mathbf{x}=\mathbf{x}^{(t)}$ as the standard minimization rule
$\gamma_{\textrm{opt}}^{(t)}$ (\ref{eq:minimization-rule-1}) does
in (\ref{eq:minimization-rule-2}), and this decreasing property is
stated in the following proposition.
\begin{prop}
Given $\tilde{\mathbf{x}}^{(t)}(\gamma)$ and $\gamma^{(t)}$ defined
in (\ref{eq:x_tilde}) and (\ref{eq:stepsize-1}), respectively, the
following holds:
\[
L^{(t)}(\tilde{\mathbf{x}}^{(t)}(\gamma^{(t)}))\leq L^{(t)}(\mathbf{x}^{(t)}),
\]
and equality is achieved if and only if $\gamma^{(t)}=0$.\end{prop}
\begin{IEEEproof}
Denote the objective function in (\ref{eq:stepsize-1}) as $\overline{L}^{(t)}(\tilde{\mathbf{x}}^{(t)}(\gamma))$.
It follows from (\ref{eq:triangle-inequality}) that \addtocounter{equation}{1}
\begin{equation}
L^{(t)}(\tilde{\mathbf{x}}^{(t)}(\gamma^{(t)}))-L^{(t)}(\mathbf{x}^{(t)})\leq\overline{L}^{(t)}(\tilde{\mathbf{x}}^{(t)}(\gamma^{(t)})),\label{eq:inequality-1}
\end{equation}
and equality in (\ref{eq:inequality-1}) is achieved when $\gamma^{(t)}=0$
and $\gamma^{(t)}=1$.

Besides, it follows from the definition of $\gamma^{(t)}$ that
\begin{equation}
\overline{L}^{(t)}(\tilde{\mathbf{x}}^{(t)}(\gamma^{(t)}))\leq\overline{L}^{(t)}(\tilde{\mathbf{x}}^{(t)}(\gamma))\bigr|_{\gamma=0}=L^{(t)}(\mathbf{x}^{(t)}).\label{eq:inequality-2}
\end{equation}
Since the optimization problem in (\ref{eq:stepsize-1}) has a unique
optimal solution $\gamma^{(t)}$ given by (\ref{eq:stepsize-2}),
equality in (\ref{eq:inequality-2}) is achieved if and only if $\gamma^{(t)}=0$.
Finally, combining (\ref{eq:inequality-1}) and (\ref{eq:inequality-2})
yields the conclusion stated in the proposition.
\end{IEEEproof}
The signaling required to perform (\ref{eq:stepsize-2}) (and also
(\ref{eq:x-hat-scalar})) will be discussed in Section \ref{sec:Implementation-and-Extensions}.

\textbf{Step 3 (Dynamic Reset): }In this step, the next estimate $\mathbf{x}^{(t+1)}$
is defined based on $\tilde{\mathbf{x}}^{(t)}(\gamma^{(t)})$ given
in (\ref{eq:x_tilde}) and (\ref{eq:stepsize-2}). We first remark
that $\tilde{\mathbf{x}}^{(t)}(\gamma^{(t)})$ is not necessarily
the solution of the optimization problem in (\ref{eq:l1-RLS}), i.e.,
\[
L^{(t)}(\tilde{\mathbf{x}}^{(t)}(\gamma^{(t)}))\geq L^{(t)}(\mathbf{x}_{\textrm{lasso}}^{(t)})=\min_{\mathbf{x}}L^{(t)}(\mathbf{x}).
\]
This is because $\mathbf{x}$ is updated only once from $\mathbf{x}=\mathbf{x}^{t}$
to $\mathbf{x}=\tilde{\mathbf{x}}^{(t)}(\gamma^{(t)})$, which in
general can be further improved unless $\tilde{\mathbf{x}}^{(t)}(\gamma^{(t)})$
already minimizes $L^{(t)}(\mathbf{x})$, i.e., $\tilde{\mathbf{x}}^{(t)}(\gamma^{(t)})=\mathbf{x}_{\textrm{lasso}}^{(t)}$.

The definitions of $L^{(t)}(\mathbf{x})$ and $\mathbf{x}_{\textrm{lasso}}^{(t)}$
in (\ref{eq:L^t(x)})-(\ref{eq:l1-RLS}) reveal that
\[
0=L^{(t)}(\mathbf{x})\bigr|_{\mathbf{x}=\mathbf{0}}\geq L^{(t)}(\mathbf{x}_{\textrm{lasso}}^{(t)}),\, t=1,2,\ldots.
\]
However, $L^{(t)}(\tilde{\mathbf{x}}^{(t)}(\gamma^{(t)}))$ may be
larger than 0 and $\tilde{\mathbf{x}}^{(t)}(\gamma^{(t)})$ is not
necessarily better than the point $\mathbf{0}$. Therefore we define
the next estimate $\mathbf{x}^{(t+1)}$ to be the best between the
two points $\tilde{\mathbf{x}}^{(t)}(\gamma^{(t)})$ and $\mathbf{0}$:
\begin{align}
\mathbf{x}^{(t+1)} & =\underset{\mathbf{x}\in\left\{ \tilde{\mathbf{x}}^{(t+1)},\mathbf{0}\right\} }{\arg\min}L^{(t)}(\mathbf{x})\nonumber \\
 & =\begin{cases}
\tilde{\mathbf{x}}^{(t)}(\gamma^{(t)}), & \textrm{ if }L^{(t)}(\tilde{\mathbf{x}}^{(t)}(\gamma^{(t)}))\leq L^{(t)}(\mathbf{0})=0,\\
\mathbf{0}, & \textrm{ otherwise, }
\end{cases}\label{eq:x^(t+1)}
\end{align}
and it is straightforward to infer the following relationship among
$\mathbf{x}^{(t)}$, $\tilde{\mathbf{x}}^{(t)}(\gamma^{(t)})$, $\mathbf{x}^{(t+1)}$
and $\mathbf{x}_{\textrm{lasso}}^{(t)}$:
\[
L^{(t)}(\mathbf{x}^{(t)})\geq L^{(t)}(\tilde{\mathbf{x}}^{(t)}(\gamma^{(t)}))\geq L^{(t)}(\mathbf{x}^{(t+1)})\geq L^{(t)}(\mathbf{x}_{\textrm{lasso}}^{(t)}).
\]
Moreover, the dynamic reset (\ref{eq:x^(t+1)}) guarantees that
\begin{equation}
\mathbf{x}^{(t+1)}\in\bigl\{\mathbf{x}:L^{(t)}(\mathbf{x})\leq0\bigr\},\; t=1,2,\ldots,\label{eq:low-level-set}
\end{equation}
Since $\lim_{t\rightarrow\infty}\mathbf{G}^{(t)}\succ\mathbf{0}$
and $\mathbf{b}^{(t)}$ converges from Assumptions (A1)-(A2), (\ref{eq:low-level-set})
guarantees that $\left\{ \mathbf{x}^{(t)}\right\} $ is a bounded
sequence.

\begin{algorithm}[t]
\textbf{Initialization:} $\mathbf{x}^{(1)}=\mathbf{0}$, $t=1$.

At each time instance $t=1,2,\ldots$:

\textbf{Step 1: }Calculate $\hat{\mathbf{x}}^{(t)}$ according to
(\ref{eq:x-hat-scalar}).

\textbf{Step 2: }Calculate $\gamma^{(t)}$ according to (\ref{eq:stepsize-2}).

\textbf{Step 3-1: }Calculate $\tilde{\mathbf{x}}^{(t)}(\gamma^{(t)})$
according to (\ref{eq:x_tilde}).

\textbf{Step 3-2: }Update $\mathbf{x}^{(t+1)}$ according to (\ref{eq:x^(t+1)}).

\protect\caption{\hspace{-0.1cm}\textbf{: }The Online Parallel Algorithm\label{alg:Iterative-Algorithm}}
\end{algorithm}

To summarize the above development, the proposed online parallel algorithm
is formally described in Algorithm \ref{alg:Iterative-Algorithm},
and its convergence properties are given in the following theorem.
\begin{thm}
[\bf{Strong Consistency}]\label{thm:convergence}Suppose Assumptions
(A1)-(A3) as well as the following assumptions are satisfied:

\begin{enumerate}

\item[\em{(A4)}] Both $\mathbf{g}_{n}$ and $v_{n}$ have bounded
moments;

\item[\em{(A5)}] $G_{kk}^{(t)}+c_{k}^{(t)}\geq c$ for some $c>0$;

\item[\em{(A6)}] The sequence $\{\mu^{(t)}\}$ is nonincreasing,
i.e., $\mu^{(t)}\geq\mu^{(t+1)}$.

\end{enumerate}Then $\mathbf{x}^{(t)}$ is a strongly consistent
estimator of $\mathbf{x}^{\star}$, i.e., $\lim_{t\rightarrow\infty}\mathbf{x}^{(t)}=\mathbf{x}^{\star}$
with probability 1. \end{thm}
\begin{IEEEproof}
See Appendix \ref{sec:Proof-of-Theorem-convergence}.
\end{IEEEproof}
Assumption (A4) is standard on random variables and is usually satisfied
in practice. We can see from Assumption (A5) that if there already
exists some $c>0$ such that $G_{kk}^{(t)}\geq c$ for all $t$, the
quadratic proximal term in (\ref{eq:x-hat-definition}) is no longer
needed, i.e., we can set $c_{k}^{(t)}=0$ without affecting convergence.
This is the case when $t$ is sufficiently large because $\lim_{t\rightarrow\infty}\mathbf{G}^{(t)}\succ\mathbf{0}$.
In practice it may be difficult to decide if $t$ is large enough,
so we can just assign a small value to $c_{k}^{(t)}$ for all $t$
in order to guarantee the convergence. As for Assumption (A6), it
is satisfied by the previously mentioned choices of $\mu^{(t)}$,
e.g., $\mu^{(t)}=\alpha/t^{\beta}$ with $\alpha>0$ and $0.5\leq\beta\leq1$.

Theorem \ref{thm:convergence} establishes that there is no loss of
strong consistency if at each time instance, (\ref{eq:l1-RLS}) is
solved only approximately by updating all elements simultaneously
based on best-response only once. In what follows, we comment on some
of the desirable features of Algorithm \ref{alg:Iterative-Algorithm}
that make it appealing in practice:

i) Algorithm \ref{alg:Iterative-Algorithm} belongs to the class of
parallel algorithms where all elements are updated simultaneously
each time. Compared with sequential algorithms where only one element
is updated at each time instance \cite{Angelosante2010}, the improvement
in convergence speed is notable, especially when the signal dimension
is large.

ii) Algorithm \ref{alg:Iterative-Algorithm} is easy to implement
and suitable for online implementation, since both the computations
of the best-response and the stepsize have closed-form expressions.
With the simplified minimization stepsize rule, notable decrease in
objective function value is achieved after each variable update, and
the trouble of tuning the decay rate of the diminishing stepsize as
required in \cite{Scutarib} is also saved. Most importantly, the
algorithm may not converge under decreasing stepsizes.

iii) Algorithm \ref{alg:Iterative-Algorithm} converges under milder
assumptions than state-of-the-art algorithms. The regression vector
$\mathbf{g}_{n}$ and noise $v_{n}$ do not need to be uniformly bounded,
which is required in \cite{Mairal2010,Razaviyayn} and which is not
satisfied in case of unbounded distribution, e.g., the Gaussian distribution.

\section{\label{sec:Implementation-and-Extensions}Implementation and Extensions}

\subsection{A special case: $\mathbf{x}^{\star}\geq\mathbf{0}$}

The proposed Algorithm \ref{alg:Iterative-Algorithm} can be further
simplified if $\mathbf{x}^{\star}$, the signal to be estimated, has
additional properties. For example, in the context of CR studied in
\cite{Kim2011a}, $\mathbf{x}^{\star}$ represents the power vector
and it is by definition always nonnegative. In this case, a nonnegative
constraint on $x_{k}$ in (\ref{eq:x-hat-definition}) is needed:
\[
\hat{x}_{k}^{(t)}=\underset{x_{k}\geq0}{\arg\min}\;\Bigl\{ L^{(t)}(x_{k},\mathbf{x}_{-k}^{(t)})+\frac{1}{2}c_{k}^{(t)}(x_{k}-x_{k}^{(t)})^{2}\Bigr\},\quad\forall\, k,
\]
 and the best-response $\hat{x}_{k}^{(t)}$ in (\ref{eq:x-hat-scalar})
is simplified to
\[
\hat{x}_{k}^{(t)}=\frac{\left[r_{k}^{(t)}+c_{k}^{(t)}x_{k}^{(t)}-\mu^{(t)}\right]^{+}}{G_{kk}^{(t)}+c_{k}^{(t)}},\; k=1,\ldots,K.
\]
Furthermore, since both $\mathbf{x}^{(t)}$ and $\hat{\mathbf{x}}^{(t)}$
are nonnegative, we have
\[
\mathbf{x}^{(t)}+\gamma(\hat{\mathbf{x}}^{(t)}-\mathbf{x}^{(t)})\geq\mathbf{0},\;0\leq\gamma\leq1,
\]
and
\[
\begin{split}\bigl\Vert\mathbf{x}^{(t)}+\gamma(\hat{\mathbf{x}}^{(t)}-\mathbf{x}^{(t)})\bigr\Vert_{1} & =\sum_{k=1}^{K}|x_{k}^{(t)}+\gamma(\hat{x}_{k}^{(t)}-x_{k}^{(t)})|\\
 & =\sum_{k=1}^{K}x_{k}^{(t)}+\gamma(\hat{x}_{k}^{(t)}-x_{k}^{(t)}).
\end{split}
\]
Therefore the \emph{standard }minimization rule (\ref{eq:minimization-rule-1})
can be adopted directly and the stepsize is accordingly given as
\[
\gamma^{(t)}=\left[-\frac{(\mathbf{G}^{(t)}\mathbf{x}^{(t)}-\mathbf{b}^{(t)}+\mu^{(t)}\mathbf{1})^{T}(\hat{\mathbf{x}}^{(t)}-\mathbf{x}^{(t)})}{(\hat{\mathbf{x}}^{(t)}-\mathbf{x}^{(t)})^{T}\mathbf{G}^{(t)}(\hat{\mathbf{x}}^{(t)}-\mathbf{x}^{(t)})}\right]_{0}^{1},
\]
where $\mathbf{1}$ is a vector with all elements equal to 1.

\subsection{Implementation issues and complexity analysis}

Algorithm \ref{alg:Iterative-Algorithm} can be implemented in both
a centralized and a distributed network architecture. To ease the
exposition, we discuss the implementation issues in the context of
WSN with a total number of $N$ sensors. The discussion for CR is
similar and thus not duplicated here.

\emph{Network with a fusion center: }The fusion center performs the
computation of (\ref{eq:x-hat-scalar}) and (\ref{eq:stepsize-2}).
To do this, the signaling from sensors to the fusion center is required:
at each time instance $t$, each sensor $n$ sends $(\mathbf{g}_{n}^{(t)},y_{n}^{(t)})\in\mathbb{R}^{K+1}$
to the fusion center. Note that $\mathbf{G}^{(t)}$ and $\mathbf{b}^{(t)}$
defined in (\ref{eq:A-and-b}) can be updated recursively:\begin{subequations}\label{eq:recursive-update}
\begin{align}
\mathbf{G}^{(t)} & =\frac{t-1}{t}\mathbf{G}^{(t-1)}+\frac{1}{t}\sum_{n=1}^{N}\mathbf{g}_{n}^{(t)}(\mathbf{g}_{n}^{(t)})^{T},\label{eq:recursive-update-1}\\
\mathbf{b}^{(t)} & =\frac{t-1}{t}\mathbf{b}^{(t-1)}+\frac{1}{t}\sum_{n=1}^{N}y_{n}^{(t)}\mathbf{g}_{n}^{(t)}.\label{eq:recursive-update-2}
\end{align}
\end{subequations}After updating $\mathbf{x}$ according to (\ref{eq:x-hat-scalar})
and (\ref{eq:stepsize-2}), the fusion center broadcasts $\mathbf{x}^{(t+1)}\in\mathbb{R}^{K}$
to all sensors.

We discuss next the computational complexity of Algorithm \ref{alg:Iterative-Algorithm}.
Note that in (\ref{eq:recursive-update}), the normalization by $t$
is not really computationally necessary because they appear in both
the numerator and denominator and thus cancel each other in the division
in (\ref{eq:x-hat-scalar}) and (\ref{eq:stepsize-2}). Computing
(\ref{eq:recursive-update-1}) requires $(N+1)(K^{2}+K)/2$ multiplications
and $(K^{2}+K)/2$ additions. To perform (\ref{eq:recursive-update-2}),
$(N+1)K$ multiplications and $NK$ additions are needed. Associated
with the computation of $\mathbf{r}^{(t)}$ in (\ref{eq:x-hat-r_t})
are $K^{2}+K$ multiplications and $2K$ additions. Then in (\ref{eq:x-hat-scalar}),
$5K$ additions, $K$ multiplications and $K$ divisions are required.
The projection in (\ref{eq:x-hat-scalar}) also requires $2K$ number
comparisons. To compute (\ref{eq:stepsize-2}), first note that $\mathbf{G}^{(t)}\mathbf{x}^{(t)}-\mathbf{b}^{(t)}$
can be recovered from $\mathbf{r}^{(t)}$ because $\mathbf{G}^{(t)}\mathbf{x}^{(t)}-\mathbf{b}^{(t)}=\mathbf{d}(\mathbf{G}^{(t)})\circ\mathbf{x}^{(t)}-\mathbf{r}^{(t)}$
and computing $\left\Vert \mathbf{x}\right\Vert _{1}$ requires $K$
number comparisons and $K-1$ additions, so what are requested in
total are $K^{2}+4K$ multiplications, $6K-2$ additions, 1 addition,
and $2K+2$ number comparisons (the projection needs at most 2 number
comparisons). The above analysis is summarized in Table \ref{tab:complexity},
and one can see that the complexity is at the order of $K^{2}$, which
is as same as traditional RLS \cite[Ch. 14]{Sayed2008Adaptive}.

\begin{table}[t]
\center%
\begin{tabular}{|c|c|c|c|c|}
\hline
equation & $+$ & $\times$ & $\div$ & min$(a,b)$\tabularnewline
\hline
(\ref{eq:recursive-update-1}) & $(K^{2}+K)/2$ & $(N+1)(K^{2}+K)/2$ & --- & ---\tabularnewline
\hline
(\ref{eq:recursive-update-2}) & $NK$ & $(N+1)K$ & --- & ---\tabularnewline
\hline
(\ref{eq:x-hat-r_t}) & $2K$ & $K^{2}+K$ & --- & ---\tabularnewline
\hline
(\ref{eq:x-hat-scalar}) & $5K$ & $K$ & $K$ & $2K$\tabularnewline
\hline
(\ref{eq:stepsize-2}) & $6K-2$ & $K^{2}+4K$ & 1 & $2K+2$\tabularnewline
\hline
\end{tabular}\protect\caption{\label{tab:complexity}Computational Complexity of Algorithm \ref{alg:Iterative-Algorithm}}

\vspace{-1cm}
\end{table}

\emph{Network without a fusion center: }In this case, the computational
tasks are evenly distributed among the sensors and the computation
of (\ref{eq:x-hat-scalar}) and (\ref{eq:stepsize-2}) is performed
locally by each sensor at the price of (limited) signaling exchange
among different sensors.

We first define the following sensor-specific variables $\mathbf{G}_{n}^{(t)}$
and $\tilde{\mathbf{b}}_{n}^{(t)}$ for sensor $n$ as follows:
\[
\mathbf{G}_{n}^{(t)}\triangleq\frac{1}{t}\sum_{\tau=1}^{t}\mathbf{g}_{n}^{(\tau)}(\mathbf{g}_{n}^{(\tau)})^{T},\textrm{ and }\mathbf{b}_{n}^{(t)}=\frac{1}{t}\sum_{\tau=1}^{t}y_{n}^{(t)}\mathbf{g}_{n}^{(t)},
\]
so that $\mathbf{G}^{(t)}=\sum_{n=1}^{N}\mathbf{G}_{n}^{(t)}$ and
$\mathbf{b}^{(t)}=\sum_{n=1}^{N}\mathbf{b}_{n}^{(t)}$. Note that
$\mathbf{G}_{n}^{(t)}$ and $\mathbf{b}_{n}^{(t)}$ can be computed
\emph{locally} by sensor $n$ and no signaling exchange is required.
It is also easy to verify that, similar to (\ref{eq:recursive-update}),
$\mathbf{G}_{n}^{(t)}$ and $\mathbf{b}_{n}^{(t)}$ can be updated
recursively by sensor $n$, so the sensors do not have to store all
past data.

The message passing among sensors in carried out in two phases. Firstly,
for sensor $n$, to perform (\ref{eq:x-hat-scalar}), $\mathbf{d}(\mathbf{G}^{(t)})$
and $\mathbf{r}^{(t)}$ are required, and they can be decomposed as
follows:\begin{subequations}
\begin{align}
\mathbf{d}(\mathbf{G}^{(t)}) & =\sum_{n=1}^{N}\mathbf{d}(\mathbf{G}_{n}^{(t)})\in\mathbb{R}^{K},\label{eq:information-1}\\
\mathbf{G}^{(t)}\mathbf{x}^{(t)}-\mathbf{b}^{(t)} & =\sum_{n=1}^{N}\left(\mathbf{G}_{n}^{(t)}\mathbf{x}^{(t)}-\mathbf{b}_{n}^{(t)}\right)\in\mathbb{R}^{K}.\label{eq:information-2}
\end{align}
 Furthermore, to determine the stepsize as in (\ref{eq:stepsize-2})
and to compare $L^{(t)}(\tilde{\mathbf{x}}^{(t)}(\gamma^{(t)}))$
with 0, the following variables are required at sensor $n$:
\begin{equation}
\mathbf{G}^{(t)}\mathbf{x}^{(t)}=\sum_{n=1}^{N}\mathbf{G}_{n}^{(t)}\mathbf{x}^{(t)}\in\mathbb{R}^{K}\label{eq:information-3}
\end{equation}
\begin{equation}
\mathbf{G}^{(t)}\hat{\mathbf{x}}^{(t)}=\sum_{n=1}^{N}\mathbf{G}_{n}^{(t)}\hat{\mathbf{x}}^{(t)}\in\mathbb{R}^{K},\label{eq:information-4}
\end{equation}
and
\begin{align}
(\mathbf{G}^{(t)}\mathbf{x}^{(t)}-\mathbf{b}^{(t)})^{T} & (\hat{\mathbf{x}}^{(t)}-\mathbf{x}^{(t)})\nonumber \\
= & \left(\sum_{n=1}^{N}(\mathbf{G}_{n}^{(t)}\mathbf{x}^{(t)}-\mathbf{b}_{n}^{(t)})\right)^{T}(\hat{\mathbf{x}}^{(t)}-\mathbf{x}^{(t)}),\label{eq:information-5}
\end{align}
\end{subequations}but computing (\ref{eq:information-5}) does not
require any additional signaling since $\sum_{n=1}^{N}(\mathbf{G}_{n}^{(t)}\mathbf{x}^{(t)}-\mathbf{b}_{n}^{(t)})$
is already available from (\ref{eq:information-2}). Note that $L^{(t)}(\tilde{\mathbf{x}}^{(t)}(\gamma^{(t)}))$
can be computed from (\ref{eq:information-2})-(\ref{eq:information-4})
because
\[
\begin{split}L^{(t)}(\tilde{\mathbf{x}}^{(t)}(\gamma^{(t)}))=\; & \frac{1}{2}(\tilde{\mathbf{x}}^{(t)}(\gamma^{(t)}))^{T}\mathbf{G}^{(t)}\tilde{\mathbf{x}}^{(t)}(\gamma^{(t)})\\
 & \quad-(\mathbf{b}^{(t)})^{T}\tilde{\mathbf{x}}^{(t)}(\gamma^{(t)})+\mu^{(t)}\bigl\Vert\tilde{\mathbf{x}}^{(t)}(\gamma^{(t)})\bigr\Vert_{1}\\
=\frac{1}{2}(\tilde{\mathbf{x}}^{(t)}( & \gamma^{(t)}))^{T}(\mathbf{G}^{(t)}\tilde{\mathbf{x}}^{(t)}(\gamma^{(t)})-2\mathbf{b}^{(t)})\\
+ & \mu^{(t)}\bigl\Vert\tilde{\mathbf{x}}^{(t)}(\gamma^{(t)})\bigr\Vert_{1}\\
=\frac{1}{2}(\tilde{\mathbf{x}}^{(t)}( & \gamma^{(t)}))^{T}(2(\mathbf{G}^{(t)}\mathbf{x}^{(t)}-\mathbf{b}^{(t)})-\mathbf{G}^{(t)}\mathbf{x}^{(t)}\\
+ & \gamma^{(t)}\mathbf{G}^{(t)}(\hat{\mathbf{x}}^{(t)}-\mathbf{x}^{(t)}))+\mu^{(t)}\bigl\Vert\tilde{\mathbf{x}}^{(t)}(\gamma^{(t)})\bigr\Vert_{1},
\end{split}
\]
where $\mathbf{G}^{(t)}\mathbf{x}^{(t)}-\mathbf{b}^{(t)}$ comes from
(\ref{eq:information-2}), $\mathbf{G}^{(t)}\mathbf{x}^{(t)}$ comes
from (\ref{eq:information-3}), and $\mathbf{G}^{(t)}(\hat{\mathbf{x}}^{(t)}-\mathbf{x}^{(t)})$
comes from (\ref{eq:information-3})-(\ref{eq:information-4}).

To summarize, in the first phase, each node needs to exchange $(\mathbf{d}(\mathbf{G}_{n}^{(t)}),\mathbf{G}_{n}^{(t)}\mathbf{x}^{(t)}-\mathbf{b}_{n}^{(t)})\in\mathbb{R}^{2K\times1}$,
while in the second phase, the sensors need to exchange $(\mathbf{G}_{n}^{(t)}\mathbf{x}^{(t)},\mathbf{G}_{n}^{(t)}\hat{\mathbf{x}}^{(t)})\in\mathbb{R}^{2K\times1}$;
thus the total signaling at each time instance is a vector of the
size $4K$. The signaling exchange can be implemented in a distributed
manner by, for example, consensus algorithms, which converge if the
graph representing the links among the sensors is connected. A detailed
discussion, however, is beyond the scope of this paper, and interested
readers are referred to \cite{Barbarossa2013} for a more comprehensive
introduction.

Now we compare Algorithm \ref{alg:Iterative-Algorithm} with state-of-the-art
distributed algorithms in terms of signaling exchange.

1) The signaling exchange of Algorithm \ref{alg:Iterative-Algorithm}
is much less than that in \cite{Barbarossa2013}. In \cite[Alg. A.5]{Barbarossa2013},
problem (\ref{eq:l1-RLS}) is solved completely for each time instance,
so it is essentially a double layer algorithm: in the inner layer,
an iterative algorithm is used to solve (\ref{eq:l1-RLS}) while in
the outer layer $t$ is increased to $t+1$ and (\ref{eq:l1-RLS})
is solved again. In each iteration of the inner layer, a vector of
the size $2K$ is exchanged among the sensors, and this is repeated
until the termination of the inner layer (which typically takes many
iterations), leading to a much heavier signaling burden than the proposed
algorithm.

2) A distributed implementation of the online sequential algorithm
in \cite{Angelosante2010} is also proposed in \cite{Kim2011a}. It
is a double layer algorithm, and in each iteration of the inner layer,
a vector with the same order of size as Algorithm \ref{alg:Iterative-Algorithm}
is exchanged among the sensors. Similar to \cite{Barbarossa2013},
this has to be repeated until the convergence of the inner layer.

We also remark that in consensus-based distributed algorithms \cite{Zeng2011,Kim2011a,Barbarossa2013},
the estimate decision of each sensor depends mainly on its own local
information and those local estimates maintained by different sensors
are usually different, which may lead to conflicting interests of
the PUs (i.e., some may correctly detect the presence of PUs but some
may not, so the PUs may still be interfered); an agreement (i.e.,
convergence) is reached only when $t$ goes to infinity \cite{boyd2004convex}.
By comparison, in the proposed algorithm, all sensors update the estimate
according to the same expression (\ref{eq:x-hat-scalar}) and (\ref{eq:stepsize-2})
based on the information jointly collected by all sensors, so they
have the same estimate of the unknown variable all the time and the
QoS of PUs are better protected.

\subsection{Time- and norm-weighted sparsity regularization}

For a given vector $\mathbf{x}$, its support $\mathcal{S}_{\mathbf{x}}$
is defined as the set of indices of nonzero elements:
\[
\mathcal{S}_{\mathbf{x}}\triangleq\{1\leq k\leq K:x_{k}\neq0\}.
\]
Suppose without loss of generality $\mathcal{S}_{\mathbf{x}^{\star}}=\{1,2,\ldots,\left\Vert \mathbf{x}^{\star}\right\Vert _{0}\}$,
where $\left\Vert \mathbf{x}\right\Vert _{0}$ is the number of nonzero
elements of $\mathbf{x}$. It is shown in \cite{Angelosante2010}
that the time-weighted sparsity regularization (\ref{eq:l1-RLS})
does not make $\mathbf{x}_{\textrm{lasso}}^{(t)}$ satisfy the so-called
``oracle properties'', which consist of support consistency, i.e.,
\[
\lim_{t\rightarrow\infty}\textrm{Prob}\left[\mathcal{S}_{\mathbf{x}_{\textrm{lasso}}^{(t)}}=\mathcal{S}_{\mathbf{x}^{\star}}\right]=1,
\]
and $\sqrt{t}$-estimation consistency, i.e.,
\[
\sqrt{t}(\mathbf{x}_{\textrm{lasso},1:\left\Vert \mathbf{x}^{\star}\right\Vert _{0}}^{(t)}-\mathbf{x}_{1:\left\Vert \mathbf{x}^{\star}\right\Vert _{0}}^{\star})\rightarrow_{d}\mathcal{N}(0,\sigma^{2}\mathbf{G}_{1:\left\Vert \mathbf{x}^{\star}\right\Vert _{0},1:\left\Vert \mathbf{x}^{\star}\right\Vert _{0}}),
\]
where $\rightarrow_{d}$ means convergence in distribution and $\mathbf{G}_{1:k,1:k}\in\mathbb{R}^{k\times k}$
is the upper left block of $\mathbf{G}$.

To make the estimation satisfy the oracle properties, it was suggested
in \cite{Angelosante2010} that a time- and norm-weighted lasso be
used, and the loss function $L^{(t)}(\mathbf{x})$ in (\ref{eq:L^t(x)})
be modified as follows:
\begin{align}
L^{(t)}(\mathbf{x}) & =\frac{1}{t}\sum_{\tau=1}^{t}\sum_{n=1}^{N}(y_{n}^{(\tau)}-(\mathbf{g}_{n}^{(\tau)})^{T}\mathbf{x})^{2}\nonumber \\
 & \quad\quad+\mu^{(t)}\sum_{k=1}^{K}\mathcal{W}_{\mu^{(t)}}(|x_{\textrm{rls},k}^{(t)}|)\cdot|x_{k}|,\label{eq:l1-RLS-weighted}
\end{align}
where:
\begin{itemize}
\item $\mathbf{x}_{\textrm{rls}}^{(t)}$ is given in (\ref{eq:rls});
\item $\lim_{t\rightarrow\infty}\mu^{(t)}=0$ and $\lim_{t\rightarrow\infty}\sqrt{t}\cdot\mu^{(t)}=\infty$,
so $\mu^{(t)}$ must decrease slower than $1/\sqrt{t}$;
\item The weight factor $\mathcal{W}_{\mu}(x)$ is defined as
\[
\mathcal{W}_{\mu}(x)\triangleq\begin{cases}
1, & \textrm{ if }x\leq\mu,\\
\frac{a\mu-x}{(a-1)\mu}, & \textrm{ if }\mu\leq x\leq a\mu,\\
0, & \textrm{ if }x\geq a\mu,
\end{cases}
\]
where $a>1$ is a given constant. Therefore, the value of the weight
function $\mu^{(t)}\mathcal{W}_{\mu^{(t)}}(|x_{\textrm{rls},k}^{(t)}|)$
in (\ref{eq:l1-RLS-weighted}) depends on the relative magnitude of
$\mu^{(t)}$ and $x_{\textrm{rls},k}^{(t)}$.
\end{itemize}
After replacing the universal sparsity regularization gain $\mu^{(t)}$
for element $x_{k}$ in (\ref{eq:x-hat-scalar}) and (\ref{eq:stepsize-2})
by $\mathcal{W}_{\mu^{(t)}}(\bigl|x_{\textrm{rls},k}^{(t)}\bigr|)$,
Algorithm \ref{alg:Iterative-Algorithm} can readily be applied to
estimate $\mathbf{x}^{\star}$ based on the time- and norm-weighted
loss function (\ref{eq:l1-RLS-weighted}) and the strong consistency
holds as well. To see this, we only need to verify the nonincreasing
property of the weight function $\mu^{(t)}\mathcal{W}_{\mu^{(t)}}(|x_{\textrm{rls},k}^{(t)}|)$.
We remark that when $t$ is sufficiently large, it is either $\mu^{(t)}\mathcal{W}_{\mu^{(t)}}(|x_{\textrm{rls},k}^{(t)}|)=0$
or $\mu^{(t)}\mathcal{W}_{\mu^{(t)}}(|x_{\textrm{rls},k}^{(t)}|)=\mu^{(t)}$.
This is because $\lim_{t\rightarrow\infty}\mathbf{x}_{\textrm{rls}}^{(t)}=\mathbf{x}^{\star}$
under the conditions of Lemma \ref{lem:connection}. If $x_{k}^{\star}>0$,
since $\lim_{t\rightarrow\infty}\mu^{(t)}=0$, we have for any arbitrarily
small $\epsilon>0$ some $t_{0}$ that $a\mu^{(t)}<x_{k}^{\star}-\epsilon$
for all $t\geq t_{0}$; the weight factor in this case is 0 for all
$t\geq t_{0}$, and the nonincreasing property is automatically satisfied.
If, on the other hand, $x_{k}^{\star}=0$, then $x_{\textrm{rls}}^{(t)}$
converges to $x_{k}^{\star}=0$ at a speed of $1/\sqrt{t}$ \cite{Greene2011}.
Since $\mu^{(t)}$ decreases slower than $1/\sqrt{t}$, we have for
some $t_{0}$ such that $x_{\textrm{rls},k}^{(t)}<\mu^{(t)}$ for
all $t\geq t_{0}$; in this case, $\mathcal{W}_{\mu^{(t)}}(x_{\textrm{rls},k}^{(t)})$
is equal to 1 and the weight factor is simply $\mu^{(t)}$ for all
$t\geq t_{0}$, which is nonincreasing.

\subsection{Recursive estimation of time-varying signals}

If the signal to be estimated is time-varying, the loss function (\ref{eq:L^t(x)})
needs to be modified in a way such that the new measurement samples
are given more weight than the old ones. Defining the so-called ``forgetting
factor'' $\beta$, where $0<\beta<1$, the new loss function is given
as follows \cite{Sayed2008Adaptive,Angelosante2010,Barbarossa2013}:

\begin{equation}
\min_{\mathbf{x}}\;\frac{1}{2t}\sum_{n=1}^{N}\sum_{\tau=1}^{t}\beta^{t-\tau}((\mathbf{g}_{n}^{(\tau)})^{T}\mathbf{x}-y_{n}^{(\tau)})^{2}+\mu^{(t)}\left\Vert \mathbf{x}\right\Vert _{1},\label{eq:time-varying-estimation}
\end{equation}
and as expected, when $\beta=1$, (\ref{eq:time-varying-estimation})
is as same as (\ref{eq:L^t(x)}). In this case, the only modification
to Algorithm \ref{alg:Iterative-Algorithm} is that $\mathbf{G}^{(t)}$
and $\mathbf{b}^{(t)}$ are updated according to the following recursive
rule:
\[
\begin{split}\mathbf{G}^{(t)} & =\frac{1}{t}\left((t-1)\beta\mathbf{G}^{(t-1)}+\sum_{n=1}^{N}\mathbf{g}_{n}^{(t)}(\mathbf{g}_{n}^{(t)})^{T}\right),\\
\mathbf{b}^{(t)} & =\frac{1}{t}\left((t-1)\beta\mathbf{b}^{(t-1)}+\sum_{n=1}^{N}y_{n}^{(t)}\mathbf{g}_{n}^{(t)}\right).
\end{split}
\]

For problem (\ref{eq:time-varying-estimation}), since the signal
to be estimated is time-varying, the convergence analysis in Theorem
\ref{thm:convergence} does not hold any more. However, simulation
results show there is little loss of optimality when optimizing (\ref{eq:time-varying-estimation})
only approximately by Algorithm \ref{alg:Iterative-Algorithm}. This
establishes the superiority of the proposed algorithm over the distributed
algorithm in \cite{Barbarossa2013} which solves (\ref{eq:time-varying-estimation})
exactly at the price of a large delay and a large signaling burden.
Besides, despite the lack of theoretical analysis, Algorithm \ref{alg:Iterative-Algorithm}
performs better than the online sequential algorithm \cite{Angelosante2010}
numerically, cf. Figure \ref{fig:SE-varying} in Section \ref{sec:Numerical-Results}.

\section{Numerical Results\label{sec:Numerical-Results}}

\begin{figure}[t]
\center

\includegraphics[scale=0.64]{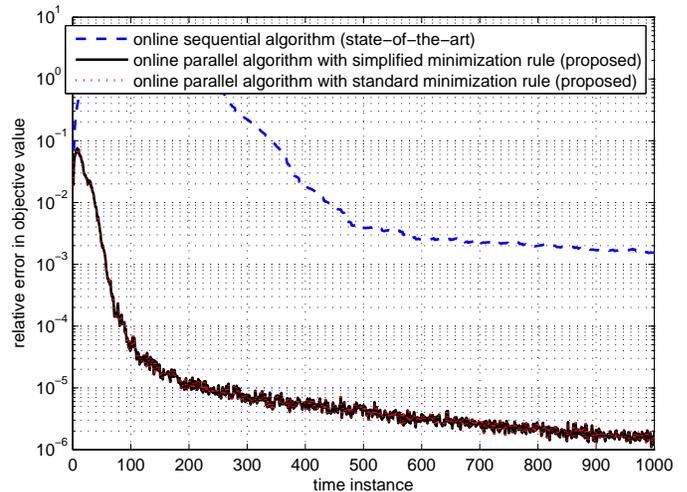}\protect\caption{\label{fig:Convergence-value}Convergence behavior in terms of objective
function value.}

\vspace{-1em}
\end{figure}

In this section, the desirable features of the proposed algorithm
are illustrated numerically.

We first test the convergence behavior of Algorithm \ref{alg:Iterative-Algorithm}
with the online sequential algorithm proposed in \cite{Angelosante2010}.
In this example, the parameters are selected as follows:
\begin{itemize}
\item $N=1$, so the subscript $n$ is omitted.
\item the dimension of $\mathbf{x}^{\star}$: $K=100$;
\item the density of $\mathbf{x}^{\star}$: 0.1;
\item Both $\mathbf{g}$ and $v$ are generated by i.i.d. standard normal
distributions: $\mathbf{g}\in\mathcal{CN}(\mathbf{0},\mathbf{I})$
and $v\in\mathcal{CN}(0,0.2)$;
\item The sparsity regularization gain $\mu^{(t)}=\sqrt{K}/t=10/t$;
\item Unless otherwise stated, the simulations results are averaged over
100 realizations.
\end{itemize}
We plot in Figure \ref{fig:Convergence-value} the relative error
in objective value $(L^{(t)}(\mathbf{x}^{(t)})-L^{(t)}(\mathbf{x}_{\textrm{lasso}}^{(t)}))/L^{(t)}(\mathbf{x}_{\textrm{lasso}}^{(t)})$
versus the time instance $t$, where 1) $\mathbf{x}_{\textrm{lasso}}^{(t)}$
is defined in (\ref{eq:l1-RLS}) and calculated by $\texttt{MOSEK}$
\cite{mosek}; 2) $\mathbf{x}^{(t)}$ is returned by Algorithm \ref{alg:Iterative-Algorithm}
in the proposed online parallel algorithm; 3) $\mathbf{x}^{(t)}$
is returned by \cite[Algorithm 1]{Angelosante2010} in online sequential
algorithm; and 4) $\mathbf{x}^{(1)}=\mathbf{0}$ for both parallel
and sequential algorithms. Note that $L^{(t)}(\mathbf{x}_{\textrm{lasso}}^{(t)})$
is by definition the lower bound of $L^{(t)}(\mathbf{x})$ and $L^{(t)}(\mathbf{x}^{(t)})-L^{(t)}(\mathbf{x}_{\textrm{lasso}}^{(t)})\geq0$
for all $t$. From Figure \ref{fig:Convergence-value} it is clear
that the proposed algorithm (black curve) converges to a precision
of $10^{-2}$ in less than 200 instances while the sequential algorithm
(blue curve) needs more than 800 instances. The improvement in convergence
speed is thus notable. If the precision is set as $10^{-4}$, the
sequential algorithm does not even converge in a reasonable number
of instances. Therefore, the proposed online parallel algorithm outperforms
in both convergence speed and solution quality.

We also evaluate in Figure \ref{fig:Convergence-value} the performance
loss incurred by the simplified minimization rule (\ref{eq:stepsize-1})
(indicated by the black curve) compared with the standard minimization
rule (\ref{eq:minimization-rule-1}) (indicated by the red curve).
It is easy to see from Figure \ref{fig:Convergence-value} that these
two curves almost coincide with each other, so the extent to which
the simplified minimization rule decreases the objective function
is nearly as same as standard minimization rule and the performance
loss is negligible.

\begin{figure}[t]
\begin{minipage}[t]{1\columnwidth}%
\begin{center}
\center\includegraphics[scale=0.65]{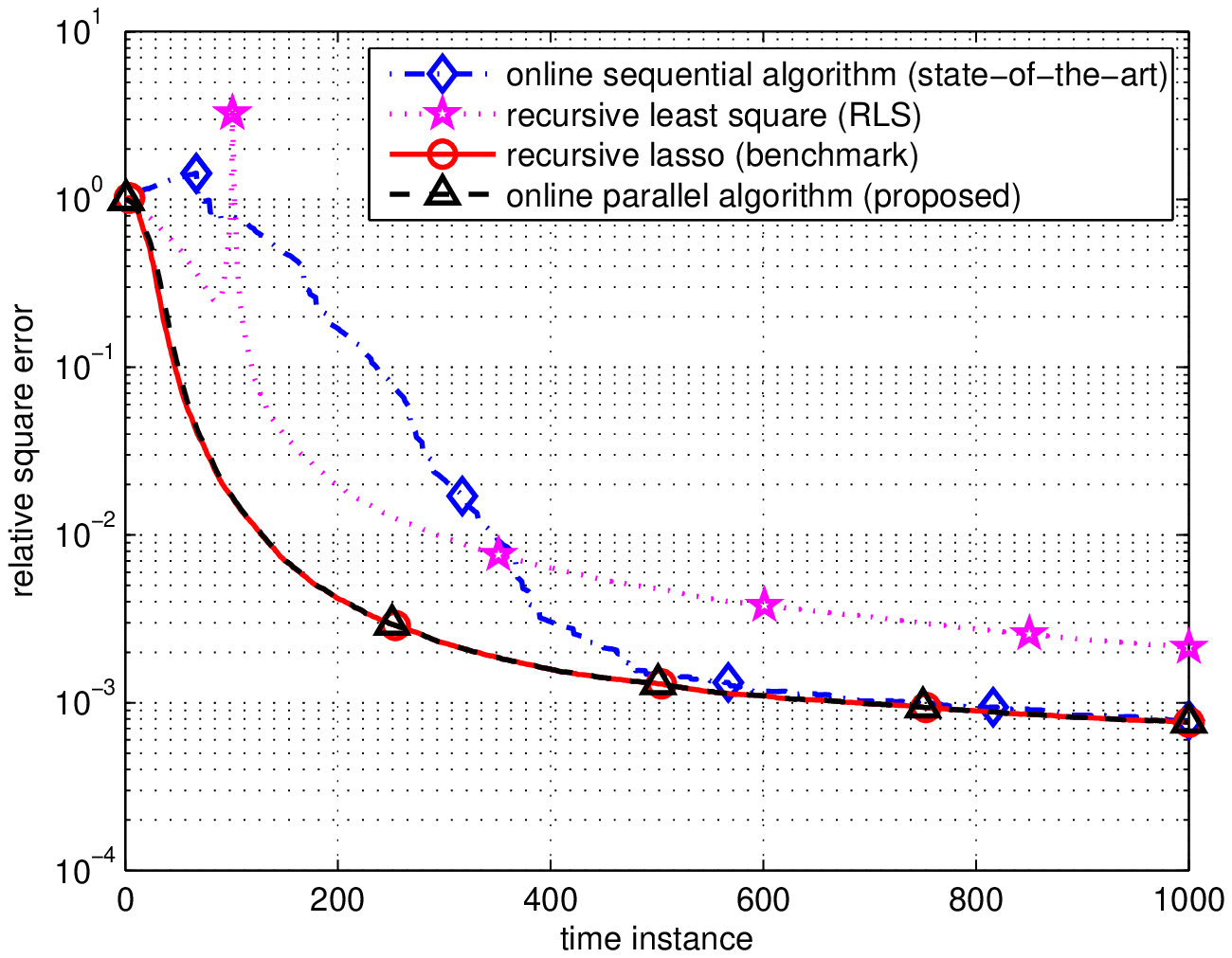}
\par\end{center}

\protect\caption{\label{fig:Convergence-MSE}Convergence behavior in terms of relative
square error.}
\end{minipage}

\vspace{0.5cm}

\begin{minipage}[t]{1\columnwidth}%
\begin{center}
\center\includegraphics[scale=0.65]{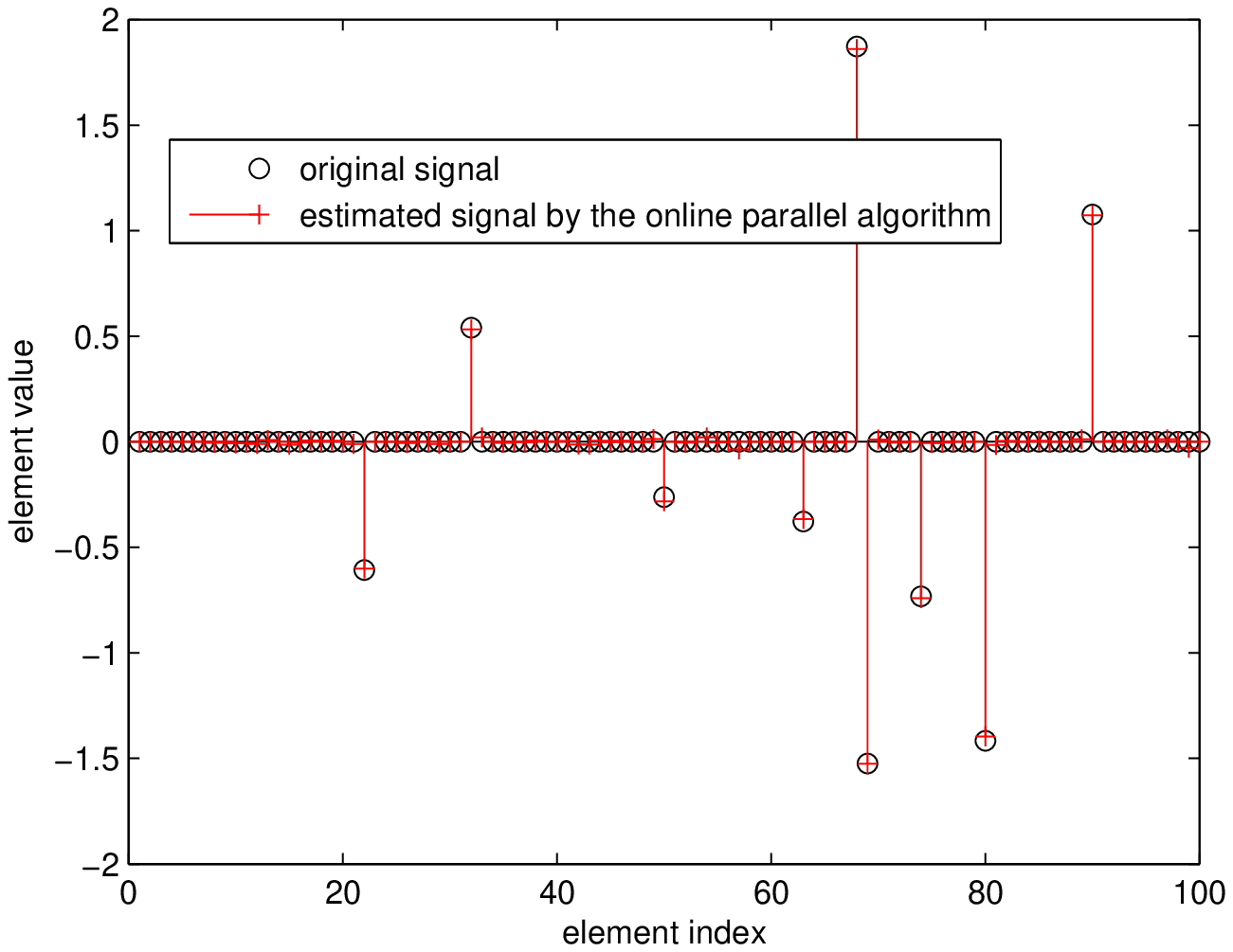}
\par\end{center}

\protect\caption{\label{fig:Comparison_signal}Comparison of original signal and estimated
signal.}
\end{minipage}
\end{figure}

Then we consider in Figure \ref{fig:Convergence-MSE} the relative
square error $\left\Vert \mathbf{x}^{(t)}-\mathbf{x}^{\star}\right\Vert _{2}^{2}/\left\Vert \mathbf{x}^{\star}\right\Vert _{2}^{2}$
versus the time instance $t$, where the benchmark is $\bigl\Vert\mathbf{x}_{\textrm{lasso}}^{(t)}-\mathbf{x}^{\star}\bigr\Vert_{2}^{2}/\left\Vert \mathbf{x}^{\star}\right\Vert _{2}^{2}$,
i.e., the recursive Lasso (\ref{eq:l1-RLS}). To compare the estimation
approaches with and without sparsity regularization, RLS in (\ref{eq:l1-RLS})
is also implemented, where a $\ell_{2}$ regularization term $(10^{-4}/t)\cdot\left\Vert \mathbf{x}\right\Vert _{2}^{2}$
is included into (\ref{eq:l1-RLS}) when $\mathbf{G}^{(t)}$ is singular.
We see that the relative square error of the proposed online parallel
algorithm quickly reaches the benchmark (recursive lasso) in about
100 instances, while the sequential algorithm needs about 800 instances.
The improvement in convergence speed is consolidated again. Another
notable feature is that, the relative square error of the proposed
algorithm is always decreasing, even in beginning instances, while
the relative square error of the sequential algorithm is not: in the
first 100 instances, the relative square error is actually increasing.
Recall the signal dimension ($K=100$), we infer that the relative
square error starts to decrease only after each element has been updated
once. What is more, estimation with sparsity regularization performs
better than the classic RLS approach because they exploit the a prior
sparsity of the to-be-estimated signal $\mathbf{x}^{\star}$. The
precision of the estimated signal by the proposed online parallel
algorithm (after 1000 time instances) is also shown element-wise in
Figure \ref{fig:Comparison_signal}, from which we observe that both
the zeros and the nonzero elements of the original signal $\mathbf{x}^{\star}$
are estimated accurately.

\begin{figure}[t]
\begin{minipage}[t]{1\columnwidth}%
\begin{center}
\center\includegraphics[scale=0.65]{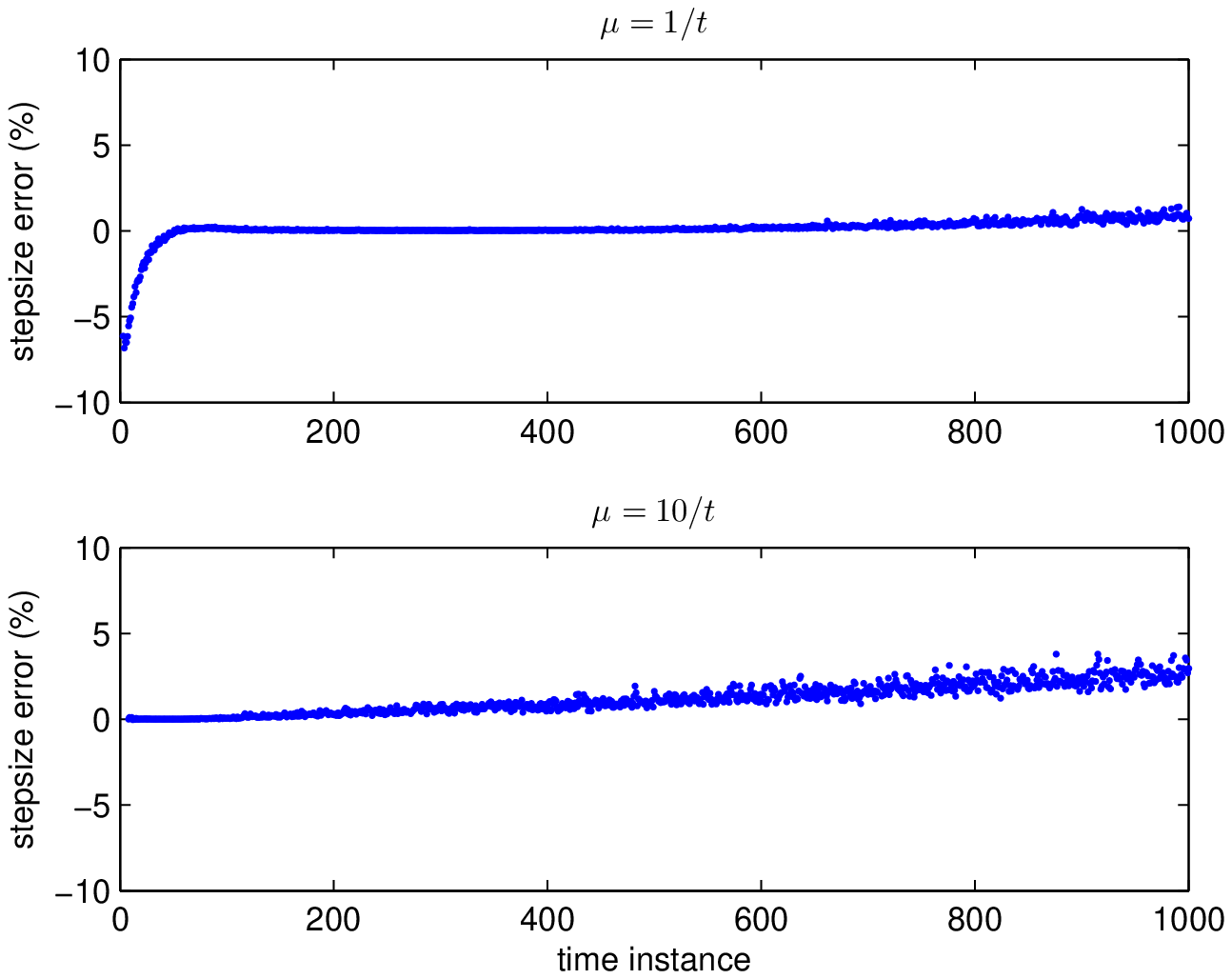}
\par\end{center}

\vspace{-0.1cm}

\protect\caption{\label{fig:Stepsize-error}Stepsize error of the simplified minimization
rule.}
\end{minipage}

\vspace{0.5cm}

\begin{minipage}[t]{1\columnwidth}%
\center\includegraphics[scale=0.66]{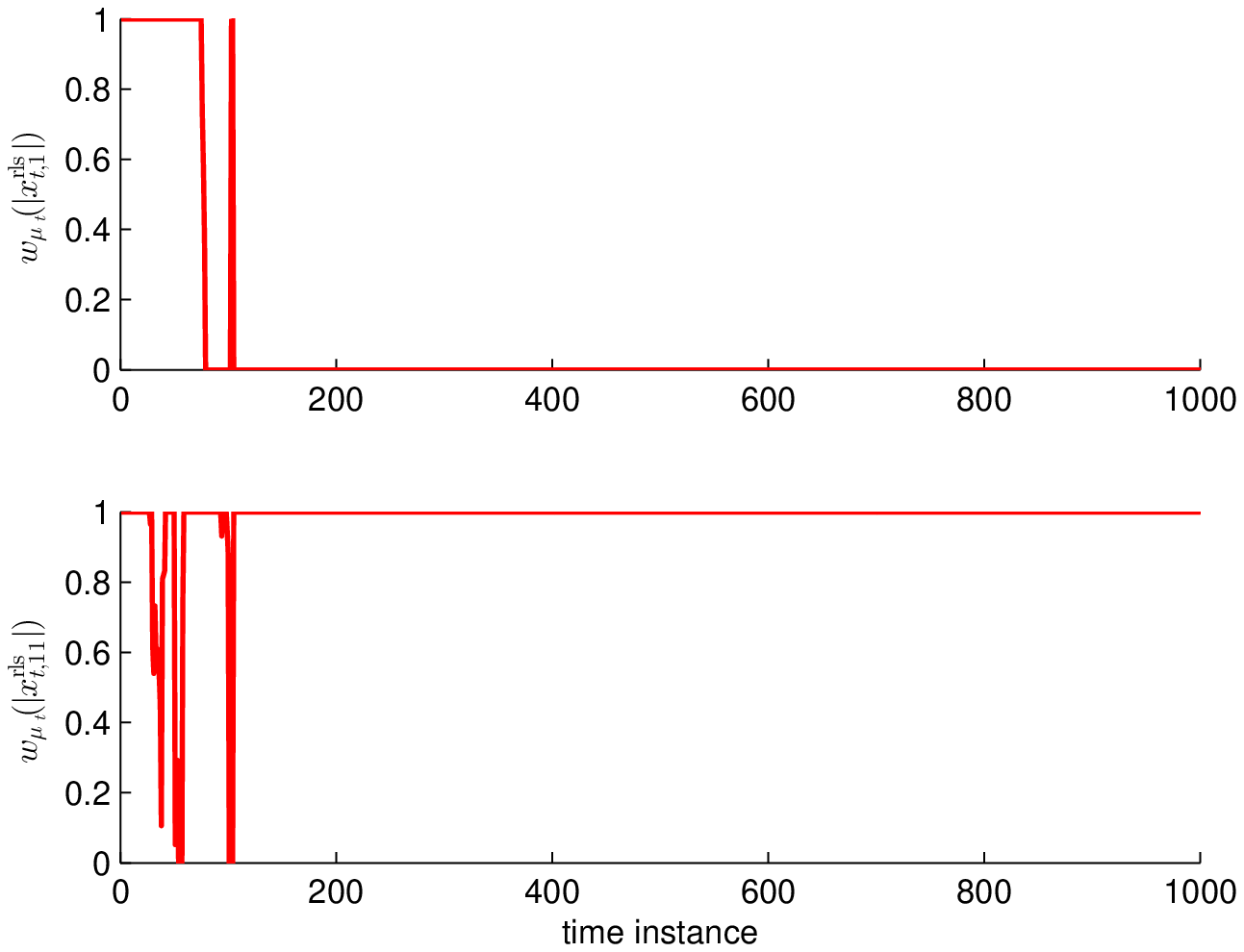}

\protect\caption{\label{fig:Weight-factor}Weight factor in time- and norm-weighted
sparsity regularization.}
\end{minipage}
\end{figure}

We now compare the proposed simplified minimization rule (cf. (\ref{eq:stepsize-1})-(\ref{eq:stepsize-2})),
coined as stepsize\_simplified, with the standard minimization rule
(cf. (\ref{eq:minimization-rule-1})), coined as stepsize\_optimal,
in terms of the stepsize error defined as follows:
\[
\frac{\textrm{stepsize\_simplified-stepsize\_optimal}}{\textrm{stepsize\_optimal}}\times100\%.
\]
In addition to the above parameter where $\mu^{(t)}=10/t$ (in the
lower subplot), we also simulate the case when $\mu^{(t)}=1/t$ (in
the upper subplot). We see from Figure \ref{fig:Stepsize-error} that
the stepsize error is reasonably small, namely, mainly in the interval
{[}-5\%,5\%{]}, while only a few of beginning instances are outside
this region, so the simplified minimization rule achieves a good trade-off
between performance and complexity. Comparing the two subplots, we
find that, as expected, the stepsize error depends on the value of
$\mu$. We can also infer that the simplified minimization rule tends
to overestimate the optimal stepsize.

\begin{figure}[t]
\center

\includegraphics[scale=0.66]{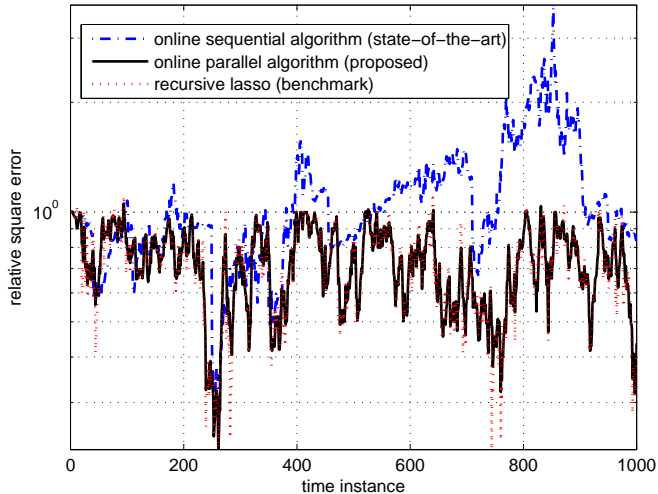}\protect\caption{\label{fig:SE-varying}Relative square error for recursive estimation
of time-varying signals}
\end{figure}

In Figure \ref{fig:Weight-factor} we simulate the weight factor $\mathcal{W}_{\mu^{(t)}}(|x_{\textrm{rls},k}^{(t)}|)$
versus the time instance $t$ in time- and norm-weighted sparsity
regularization, where $k=1$ in the upper plot and $k=11$ in the
lower plot. The parameters are as same as the above examples, except
that $\mu^{(t)}=1/t^{0.4}$ and $\mathbf{x}^{\star}$ are generated
such that the first $0.1\times K$ elements (where 0.1 is the density
of $\mathbf{x}^{\star}$) are nonzero while all other elements are
zero. The weight factors of other elements are omitted because they
exhibit similar behavior as the ones plotted in Figure \ref{fig:Weight-factor}.
As analyzed, $\mathcal{W}_{\mu^{(t)}}(|w_{\textrm{rls},1}^{(t)}|)$,
the weight factor of the first element, where $x_{1}^{\star}\neq0$,
quickly converges to 0, while $\mathcal{W}_{\mu^{(t)}}(|w_{\textrm{rls},11}^{(t)}|)$,
the weight factor of the eleventh element, where $x_{1}^{\star}=0$
, quickly converges to $1$, making the overall weight factor monotonically
decreasing, cf. (\ref{eq:l1-RLS-weighted}). Therefore the proposed
algorithm can readily be applied to the recursive estimation of sparse
signals with time- and norm-weighted regularization.

When the signal to be estimated is time-varying, the theoretical analysis
of the proposed algorithm is not valid anymore, but we can test numerically
how the proposed algorithm performs compared with the online sequential
algorithm. The time-varying unknown signal is denoted as $\mathbf{x}_{t}^{\star}$,
and it is changing according to the following law:
\[
x_{t+1,k}^{\star}=\alpha x_{t,k}^{\star}+w_{t,k},
\]
where $w_{t,k}\sim\mathcal{CN}(0,1-\alpha^{2})$ for any $k$ such
that $x_{t,k}^{\star}\neq0$, with $\alpha=0.99$ and $\beta=0.9$.
In Figure \ref{fig:SE-varying}, the relative square error $\left\Vert \mathbf{x}_{t}-\mathbf{x}_{t}^{\star}\right\Vert _{2}^{2}/\left\Vert \mathbf{x}_{t}^{\star}\right\Vert _{2}^{2}$
is plotted versus the time instance. Despite the lack of theoretical
consolidation, we observe the online parallel algorithm is almost
as same as the pseudo-online algorithm, so the inexact optimization
is not an impeding factor for the estimation accuracy. This also consolidates
the superiority of the proposed algorithm over \cite{Barbarossa2013}
where a distributed iterative algorithm is employed to solve (\ref{eq:time-varying-estimation})
exactly, which inevitably incurs a large delay and extensive signaling.

\section{Concluding Remarks\label{sec:Concluding-Remarks}}

In this paper, we have considered the recursive estimation of sparse
signals and proposed an online parallel algorithm with provable convergence.
The algorithm is based on inexact optimization with an increasing
accuracy and parallel update of all elements at each time instance,
where both the update direction and the stepsize can be calculated
in closed-form expressions. The proposed simplified minimization stepsize
rule is well motivated and easily implementable, achieving a good
trade-off between complexity and convergence speed, and avoiding the
common drawbacks of the decreasing stepsizes used in literature. Simulation
results consolidate the notable improvement in convergence speed over
state-of-the-art techniques, and they also show that the loss in convergence
speed compared with the full version (where the lasso problem is solved
exactly at each time instance) is negligible. We have also considered
numerically the recursive estimation of time-varying signals where
theoretical convergence do not necessarily hold, and the proposed
algorithm works better than state-of-the-art algorithms.

\appendices{}

\section{\label{sec:Proof-of-Theorem-convergence}Proof of Theorem \ref{thm:convergence}}
\begin{IEEEproof}
It is easy to see that $L^{(t)}$ can be divided into a smooth part
$f^{(t)}(\mathbf{x})$ and a nonsmooth part $h^{(t)}(\mathbf{x})$:\begin{subequations}\label{rem:decompose-of-L_t(x)}
\begin{align}
f^{(t)}(\mathbf{x}) & \triangleq\frac{1}{2}\mathbf{x}^{T}\mathbf{G}^{(t)}\mathbf{x}-(\mathbf{b}^{(t)})^{T}\mathbf{x},\label{eq:f(x)}\\
h^{(t)}(\mathbf{x}) & \triangleq\mu^{(t)}\left\Vert \mathbf{x}\right\Vert _{1}.\label{eq:h(x)}
\end{align}
\end{subequations}We also use $f_{k}^{(t)}(x;\mathbf{x}^{(t)})$
to denote the smooth part of the objective function in (\ref{eq:x-hat-scalar}):
\begin{equation}
f_{k}^{(t)}(x;\mathbf{x}^{(t)})\triangleq\frac{1}{2}G_{kk}^{(t)}x^{2}-r_{k}^{(t)}\cdot x+\frac{1}{2}c_{k}^{(t)}(x-x_{k}^{(t)})^{2}.\label{eq:f_t-k}
\end{equation}
Functions $f_{k}^{(t)}(x;\mathbf{x}^{(t)})$ and $f^{(t)}(\mathbf{x})$
are related according to the following equation:
\begin{equation}
f_{k}^{(t)}(x_{k};\mathbf{x}^{(t)})=f^{(t)}(x_{k},\mathbf{x}_{-k}^{(t)})+\frac{1}{2}c^{(t)}(x_{k}-x_{k}^{(t)})^{2},\label{eq:functions-relation}
\end{equation}
from which it is easy to infer that $\nabla f_{k}^{(t)}(x_{k}^{(t)};\mathbf{x}^{(t)})=\nabla_{k}f^{(t)}(\mathbf{x}^{(t)})$.
Then from the first-order optimality condition, $h^{(t)}(x_{k})$
has a subgradient $\xi_{k}^{(t)}\in\partial h^{(t)}(\hat{x}_{k}^{(t)})$
at $x_{k}=\hat{x}_{k}^{(t)}$ such that for any $x_{k}$:
\begin{equation}
(x_{k}-\hat{x}_{k}^{(t)})(\nabla f_{k}^{(t)}(\hat{x}_{k}^{(t)};\mathbf{x}^{(t)})+\xi_{k}^{(t)})\geq0,\,\forall k.\label{eq:minimum-principle}
\end{equation}

Now consider the following equation:
\begin{equation}
\begin{array}{l}
L^{(t)}(\mathbf{x}^{(t+1)})-L^{(t-1)}(\mathbf{x}^{(t)})=\smallskip\\
\qquad L^{(t)}(\mathbf{x}^{(t+1)})-L^{(t)}(\mathbf{x}^{(t)})+L^{(t)}(\mathbf{x}^{(t)})-L^{(t-1)}(\mathbf{x}^{(t)}).
\end{array}\label{eq:basic-relation}
\end{equation}
The rest of the proof consists of three parts. Firstly we prove that
there exists a constant $\eta>0$ such that $L^{(t)}(\mathbf{x}^{(t+1)})-L^{(t)}(\mathbf{x}^{(t)})\leq-\eta\left\Vert \hat{\mathbf{x}}^{(t)}-\mathbf{x}^{(t)}\right\Vert _{2}^{2}$.
Then we show that the sequence $\left\{ L^{(t)}(\mathbf{x}^{(t+1)})\right\} _{t}$
converges. Finally we prove that any limit point of the sequence $\left\{ \mathbf{x}^{(t)}\right\} _{t}$
is a solution of (\ref{eq:mmse}).

\textbf{Part 1)} Since $c_{\min}^{(t)}\geq c>0$ for all $t$ ($c_{\textrm{min}}^{(t)}$
is defined in Proposition \ref{prop:descent-direction}) from Assumption
(A5), it is easy to see from (\ref{eq:descent}) that the following
is true:
\begin{align*}
 & L^{(t)}(\mathbf{x}^{(t)}+\gamma(\hat{\mathbf{x}}^{(t)}-\mathbf{x}^{(t)}))-L^{(t)}(\mathbf{x}^{(t)})\\
\leq\; & -\gamma\Bigl(c-\frac{1}{2}\lambda_{\max}(\mathbf{G}^{(t)})\gamma\Bigr)\bigl\Vert\hat{\mathbf{x}}^{(t)}-\mathbf{x}^{(t)}\bigr\Vert_{2}^{2},\,0\leq\gamma\leq1.
\end{align*}
Since $\lambda_{\max}(\bullet)$ is a continuous function \cite{Matrix_Analysis}
and $\mathbf{G}^{(t)}$ converges to a positive definite matrix by
Assumption (A1), there exists a $\bar{\lambda}<+\infty$ such that
$\bar{\lambda}\geq\lambda_{\max}(\mathbf{G}^{(t)})$ for all $t$.
We thus conclude from the preceding inequality that for all $0\leq\lambda\leq1$:
\begin{align}
 & L^{(t)}(\mathbf{x}^{(t)}+\gamma(\hat{\mathbf{x}}^{(t)}-\mathbf{x}^{(t)}))-L^{(t)}(\mathbf{x}^{(t)})\nonumber \\
\leq & -\gamma\left(c-\frac{1}{2}\bar{\lambda}\gamma\right)\bigl\Vert\hat{\mathbf{x}}^{(t)}-\mathbf{x}^{(t)}\bigr\Vert_{2}^{2}.\label{eq:descent-lemma-4}
\end{align}

It follows from (\ref{eq:triangle-inequality}), (\ref{eq:stepsize-1})
and (\ref{eq:descent-lemma-4}) that
\begin{align}
 & L^{(t)}(\tilde{\mathbf{x}}^{(t+1)})\nonumber \\
\leq\; & f^{(t)}(\mathbf{x}^{(t)}+\gamma^{(t)}(\hat{\mathbf{x}}^{(t)}-\mathbf{x}^{(t)}))\nonumber \\
 & +(1-\gamma^{(t)})h^{(t)}(\mathbf{x}^{(t)})+\gamma^{(t)}h^{(t)}(\hat{\mathbf{x}}^{(t)})\\
\leq\; & f^{(t)}(\mathbf{x}^{(t)}+\gamma(\hat{\mathbf{x}}^{(t)}-\mathbf{x}^{(t)}))\nonumber \\
 & +(1-\gamma)h^{(t)}(\mathbf{x}^{(t)})+\gamma h^{(t)}(\hat{\mathbf{x}}^{(t)})\\
\leq\; & L^{(t)}(\mathbf{x}^{(t)})-\gamma(c-\frac{1}{2}\bar{\lambda}\gamma)\bigl\Vert\hat{\mathbf{x}}^{(t)}-\mathbf{x}^{(t)}\bigr\Vert_{2}^{2}.\label{eq:descent-lemma}
\end{align}
Since the inequalities in (\ref{eq:descent-lemma}) are true for any
$0\leq\gamma\leq1$, we set $\gamma=\min(c/\bar{\lambda},1)$. Then
it is possible to show that there is a constant $\eta>0$ such that
\begin{align}
L^{(t)}(\mathbf{x}^{(t+1)})-L^{(t)}(\mathbf{x}^{(t)}) & \leq L^{(t)}(\tilde{\mathbf{x}}^{(t+1)})-L^{(t)}(\mathbf{x}^{(t)})\nonumber \\
 & \leq-\eta\bigl\Vert\hat{\mathbf{x}}^{(t)}-\mathbf{x}^{(t)}\bigr\Vert_{2}^{2}.\label{eq:descent-property-1}
\end{align}
Besides, because of Step 3 in Algorithm \ref{alg:Iterative-Algorithm},
$\mathbf{x}^{(t+1)}$ is in the following lower level set of $L^{(t)}(\mathbf{x})$:
\begin{equation}
\mathcal{L}_{\leq0}^{(t)}\triangleq\{\mathbf{x}:L^{(t)}(\mathbf{x})\leq0\}.\label{eq:lower-level-set}
\end{equation}
Because $\left\Vert \mathbf{x}\right\Vert _{1}\geq0$ for any $\mathbf{x}$,
(\ref{eq:lower-level-set}) is a subset of
\[
\left\{ \mathbf{x}:\frac{1}{2}\mathbf{x}^{T}\mathbf{G}^{(t)}\mathbf{x}-(\mathbf{b}^{(t)})^{T}\mathbf{x}\leq0\right\} ,
\]
which is a subset of
\begin{equation}
\bar{\mathcal{L}}_{\leq0}^{(t)}\triangleq\left\{ \mathbf{x}:\frac{1}{2}\lambda_{\max}(\mathbf{G}^{(t)})\left\Vert \mathbf{x}\right\Vert _{2}^{2}-(\mathbf{b}^{(t)})^{T}\mathbf{x}\leq0\right\} .\label{eq:lower-level-set-1}
\end{equation}
Since $\mathbf{G}^{(t)}$ and $\mathbf{b}^{(t)}$ converges and $\lim_{t\rightarrow\infty}\mathbf{G}^{(t)}\succ\mathbf{0}$,
there exists a bounded set, denoted as $\mathcal{L}_{\leq0}$, such
that $\mathcal{L}_{\leq0}^{(t)}\subseteq\bar{\mathcal{L}}_{\leq0}^{(t)}\subseteq\mathcal{L}_{\leq0}$
for all $t$; thus the sequence $\{\mathbf{x}^{(t)}\}$ is bounded
and we denote its upper bound as $\bar{\mathbf{x}}$.

\textbf{Part 2)} Combining (\ref{eq:basic-relation}) and (\ref{eq:descent-property-1}),
we have the following:
\begin{align}
 & L^{(t+1)}(\mathbf{x}^{(t+2)})\negmedspace-\negmedspace L^{(t)}(\mathbf{x}^{(t+1)})\nonumber \\
\leq\, & L^{(t+1)}(\mathbf{x}^{(t+1)})\negmedspace-\negmedspace L^{(t)}(\mathbf{x}^{(t+1)})\nonumber \\
=\, & f^{(t+1)}(\mathbf{x}^{(t+1)})\negmedspace-\negmedspace f^{(t)}(\mathbf{x}^{(t+1)})\negmedspace+\negmedspace h^{(t+1)}(\mathbf{x}^{(t+1)})\negmedspace-\negmedspace h^{(t)}(\mathbf{x}^{(t+1)})\nonumber \\
\leq\, & f^{(t+1)}(\mathbf{x}^{(t+1)})\negmedspace-\negmedspace f^{(t)}(\mathbf{x}^{(t+1)}),\label{eq:step2-0}
\end{align}
where the last inequality comes from the decreasing property of $\mu^{(t)}$
by Assumption (A6). Recalling the definition of $f^{(t)}(\mathbf{x})$
in (\ref{rem:decompose-of-L_t(x)}), it is easy to see that
\begin{align*}
(t+1)(f^{(t+1)}(\mathbf{x}^{(t+1)})-f^{(t)}(\mathbf{x}^{(t+1)}))\\
=l^{(t+1)}(\mathbf{x}^{(t+1)}) & -\frac{1}{t}\sum_{\tau=1}^{t}l^{(\tau)}(\mathbf{x}^{(t+1)}),
\end{align*}
where
\[
l^{(t)}(\mathbf{x})\triangleq\sum_{n=1}^{N}(y_{n}^{(t)}-(\mathbf{g}_{n}^{(t)})^{T}\mathbf{x})^{2}.
\]

Taking the expectation of the preceding equation with respect to $\{y_{n}^{(t+1)},\mathbf{g}_{n}^{(t+1)}\}_{n=1}^{N}$,
conditioned on the natural history up to time $t+1$, denoted as $\mathcal{F}^{(t+1)}$:
\[
\begin{array}{l}
\mathcal{F}^{(t+1)}=\smallskip\\
\quad\left\{ \mathbf{x}^{(0)},\ldots,\mathbf{x}^{(t+1)},\bigl\{\mathbf{g}_{n}^{(0)},\ldots,\mathbf{g}_{n}^{(t)}\bigr\}_{n},\bigl\{ y_{n}^{(0)},\ldots,y_{n}^{(t)}\bigr\}_{n}\right\} ,
\end{array}
\]
we have
\begin{align}
 & \mathbb{E}\left[(t+1)(f^{(t+1)}(\mathbf{x}^{(t+1)})-f^{(t)}(\mathbf{x}^{(t+1)}))|\mathcal{F}^{(t+1)}\right]\nonumber \\
=\; & \mathbb{E}\left[l^{(t+1)}(\mathbf{x}^{(t+1)})|\mathcal{F}^{(t+1)}\right]-\frac{1}{t}\sum_{\tau=1}^{t}\mathbb{E}\left[l^{(\tau)}(\mathbf{x}^{(t+1)})|\mathcal{F}^{(t+1)}\right]\nonumber \\
=\; & \mathbb{E}\left[l^{(t+1)}(\mathbf{x}^{(t+1)})|\mathcal{F}^{(t+1)}\right]-\frac{1}{t}\sum_{\tau=1}^{t}l^{(\tau)}(\mathbf{x}^{(t+1)}),\label{eq:expectation-conditional}
\end{align}
where the second equality comes from the observation that $l^{(\tau)}(\mathbf{x}^{(t+1)})$
is deterministic as long as $\mathcal{F}^{(t+1)}$ is given. This
together with (\ref{eq:step2-0}) indicates that
\begin{align*}
 & \mathbb{E}\left[L^{(t+1)}(\mathbf{x}^{(t+2)})-L^{(t)}(\mathbf{x}^{(t+1)})|\mathcal{F}^{(t+1)}\right]\\
\leq\; & \mathbb{E}\left[f^{(t+1)}(\mathbf{x}^{(t+1)})-f^{(t)}(\mathbf{x}^{(t+1)})|\mathcal{F}^{(t+1)}\right]\\
\leq\; & \frac{1}{t+1}\left(\mathbb{E}\left[l^{(t+1)}(\mathbf{x}^{(t+1)})|\mathcal{F}^{(t+1)}\right]-\frac{1}{t}\sum_{\tau=1}^{t}l^{(\tau)}(\mathbf{x}^{(t+1)})\right)\\
\leq\; & \frac{1}{t+1}\left|\mathbb{E}\left[l^{(t+1)}(\mathbf{x}^{(t+1)})|\mathcal{F}^{(t+1)}\right]-\frac{1}{t}\sum_{\tau=1}^{t}l^{(\tau)}(\mathbf{x}^{(t+1)})\right|,
\end{align*}
and
\begin{align}
 & \left[\mathbb{E}\left[L^{(t+1)}(\mathbf{x}^{(t+2)})-L^{(t)}(\mathbf{x}^{(t+1)})|\mathcal{F}^{(t+1)}\right]\right]_{0}\nonumber \\
\leq\; & \frac{1}{t+1}\left|\mathbb{E}\left[l^{(t+1)}(\mathbf{x}^{(t+1)})|\mathcal{F}^{(t+1)}\right]-\frac{1}{t}\sum_{\tau=1}^{t}l^{(\tau)}(\mathbf{x}^{(t+1)})\right|\nonumber \\
\leq\; & \frac{1}{t+1}\sup_{\mathbf{x}\in\mathcal{X}}\left|\mathbb{E}\left[l^{(t+1)}(\mathbf{x})|\mathcal{F}^{(t+1)}\right]-\frac{1}{t}\sum_{\tau=1}^{t}l^{(\tau)}(\mathbf{x})\right|,\label{eq:step2-1-1}
\end{align}
where $[x]_{0}=\max(x,0)$, and $\mathcal{X}$ in (\ref{eq:step2-1-1})
with $\mathcal{X}\triangleq\{\mathbf{x}^{(1)},\mathbf{x}^{(2)},\ldots,\}$
is the complete path of $\mathbf{x}$.

Now we derive an upper bound on the expected value of the right hand
side of (\ref{eq:step2-1-1}):
\begin{align}
 & \mathbb{E}\left[\sup_{\mathbf{x}\in\mathcal{X}}\left|\mathbb{E}\left[l^{(t+1)}(\mathbf{x})|\mathcal{F}^{(t+1)}\right]-\frac{1}{t}\sum_{\tau=1}^{t}l^{(\tau)}(\mathbf{x})\right|\right]\nonumber \\
=\; & \mathbb{E}\left[\sup_{\mathbf{x}\in\mathcal{X}}\bigl|\breve{y}^{(t)}-(\mathbf{r}_{2}^{(t)})^{T}\mathbf{x}+\mathbf{x}^{T}\mathbf{R}_{3}^{(t)}\mathbf{x}\bigr|\right]\nonumber \\
\leq\; & \mathbb{E}\left[\sup_{\mathbf{x}\in\mathcal{X}}\bigl|\breve{y}^{(t)}\bigr|+\sup_{\mathbf{x}\in\mathcal{X}}\bigl|(\breve{\mathbf{b}}^{(t)})^{T}\mathbf{x}\bigr|+\sup_{\mathbf{x}\in\mathcal{X}}\bigl|\mathbf{x}^{T}\breve{\mathbf{G}}^{(t)}\mathbf{x}\bigr|\right]\nonumber \\
=\; & \mathbb{E}\left[\sup_{\mathbf{x}\in\mathcal{X}}\bigl|\breve{y}^{(t)}\bigr|\right]+\mathbb{E}\left[\sup_{\mathbf{x}\in\mathcal{X}}\bigl|(\breve{\mathbf{b}}^{(t)})^{T}\mathbf{x}\bigr|\right]+\mathbb{E}\left[\sup_{\mathbf{x}\in\mathcal{X}}\bigl|\mathbf{x}^{T}\breve{\mathbf{G}}^{(t)}\mathbf{x}\bigr|\right],\label{eq:bound-1}
\end{align}
where
\begin{align*}
\breve{y}^{(t)} & \triangleq\frac{1}{t}\sum_{\tau=1}^{t}\sum_{n=1}^{N}\left(\mathbb{E}_{y_{n}}\left[y_{n}^{2}\right]-(y_{n}^{(\tau)})^{2}\right),\\
\breve{\mathbf{b}}^{(t)} & \triangleq\frac{1}{t}\sum_{\tau=1}^{t}\sum_{n=1}^{N}2\left(\mathbb{E}_{\{y_{n},\mathbf{g}_{n}\}}\left[y_{n}\mathbf{g}_{n}\right]-y_{n}^{(\tau)}\mathbf{g}_{n}^{(\tau)}\right),\\
\breve{\mathbf{G}}^{(t)} & \triangleq\frac{1}{t}\sum_{\tau=1}^{t}\sum_{n=1}^{N}\left(\mathbb{E}_{\mathbf{g}_{n}}\left[\mathbf{g}_{n}\mathbf{g}_{n}\right]-\mathbf{g}_{n}^{(t)}\mathbf{g}_{n}^{(\tau)T}\right).
\end{align*}
Then we bound each term in (\ref{eq:bound-1}) individually. For the
first term, since $\breve{y}^{(t)}$ is independent of $\mathbf{x}^{(t)}$,
\begin{align}
\mathbb{E}\left[\sup_{\mathbf{x}\in\mathcal{X}}\bigl|\breve{y}^{(t)}\bigr|\right]=\mathbb{E}\left[\bigl|\breve{y}^{(t)}\bigr|\right] & =\mathbb{E}\left[\sqrt{(\breve{y}^{(t)})^{2}}\right]\nonumber \\
 & \leq\sqrt{\mathbb{E}\left[(\breve{y}^{(t)})^{2}\right]}\leq\sqrt{\frac{\sigma_{1}^{2}}{t}}\label{eq:bound-2}
\end{align}
for some $\sigma_{1}<\infty$, where the second equality comes from
Jensen's inequality. Because of Assumptions (A1), (A2) and (A4), $\breve{y}^{(t)}$
has bounded moments and the existence of $\sigma_{1}$ is then justified
by the central limit theorem \cite{Durrett_Probability}.

For the second term of (\ref{eq:bound-1}), we have
\begin{align*}
\mathbb{E}\left[\sup_{\mathbf{x}}\bigl|(\breve{\mathbf{b}}^{(t)})^{T}\mathbf{x}\bigr|\right] & \negmedspace\leq\negthinspace\mathbb{E}\left[\sup_{\mathbf{x}}(\bigl|\breve{\mathbf{b}}^{(t)}\bigr|)^{T}\left|\mathbf{x}\right|\right]\negmedspace\leq\negmedspace\left(\mathbb{E}\left[\bigl|\breve{\mathbf{b}}^{(t)}\bigr|\right]\right)^{\negmedspace T}\negmedspace\left|\bar{\mathbf{x}}\right|.
\end{align*}
Similar to the line of analysis of (\ref{eq:bound-2}), there exists
a $\sigma_{2}<\infty$ such that
\begin{equation}
\mathbb{E}\left[\sup_{\mathbf{x}}\bigl|(\breve{\mathbf{b}}^{(t)})^{T}\mathbf{x}\bigr|\right]\leq\left(\mathbb{E}\left[\bigl|\breve{\mathbf{b}}^{(t)}\bigr|\right]\right)^{T}\left|\bar{\mathbf{x}}\right|\leq\sqrt{\frac{\sigma_{2}^{2}}{t}}.\label{eq:bound-3}
\end{equation}

For the third term of (\ref{eq:bound-1}), we have
\begin{align}
 & \mathbb{E}\left[\sup_{\mathbf{x}\in\mathcal{X}}\bigl|\mathbf{x}^{T}\breve{\mathbf{G}}^{(t)}\mathbf{x}\bigr|\right]\nonumber \\
=\; & \mathbb{E}\left[\max_{1\leq k\leq K}\bigl|\lambda_{k}(\breve{\mathbf{G}}^{(t)})\bigr|\cdot\left\Vert \bar{\mathbf{x}}\right\Vert _{2}^{2}\right]\nonumber \\
=\; & \left\Vert \bar{\mathbf{x}}\right\Vert _{2}^{2}\cdot\mathbb{E}\left[\sqrt{\max\{\lambda_{\max}^{2}(\breve{\mathbf{G}}^{(t)}),\lambda_{\min}^{2}(\breve{\mathbf{G}}^{(t)})\}}\right]\nonumber \\
\leq\; & \left\Vert \bar{\mathbf{x}}\right\Vert _{2}^{2}\cdot\sqrt{\mathbb{E}\left[\max\{\lambda_{\max}^{2}(\breve{\mathbf{G}}^{(t)}),\lambda_{\min}^{2}(\breve{\mathbf{G}}^{(t)})\}\right]}\nonumber \\
\leq\; & \left\Vert \bar{\mathbf{x}}\right\Vert _{2}^{2}\cdot\sqrt{\mathbb{E}\left[\sum_{k=1}^{K}\lambda_{k}^{2}(\breve{\mathbf{G}}^{(t)})\right]}\nonumber \\
=\; & \left\Vert \bar{\mathbf{x}}\right\Vert _{2}^{2}\cdot\sqrt{\mathbb{E}\left[\textrm{tr}\left(\breve{\mathbf{G}}^{(t)}(\breve{\mathbf{G}}^{(t)})^{T}\right)\right]}\leq\sqrt{\frac{\sigma_{3}^{2}}{t}}\label{eq:bound-4}
\end{align}
for some $\sigma_{3}<\infty$, where the first equality comes from
the observation that $\mathbf{x}$ should align with the eigenvector
associated with the eigenvalue with largest absolute value. Then combing
(\ref{eq:bound-2})-(\ref{eq:bound-4}), we can claim that there exists
$\sigma\triangleq\sqrt{\sigma_{1}^{2}}+\sqrt{\sigma_{2}^{2}}+\sqrt{\sigma_{3}^{2}}>0$
such that
\[
\mathbb{E}\left[\sup_{\mathbf{x}\in\mathcal{X}}\left|\mathbb{E}\left[l^{(t+1)}(\mathbf{x})|\mathcal{F}^{(t+1)}\right]-\frac{1}{t}\sum_{\tau=1}^{t}l^{(\tau)}(\mathbf{x})\right|\right]\leq\frac{\sigma}{\sqrt{t}}.
\]
In view of (\ref{eq:step2-1-1}), we have
\begin{equation}
\mathbb{E}\left[\left[\mathbb{E}\left[L^{(t+1)}(\mathbf{x}^{(t+2)})-L^{(t)}(\mathbf{x}^{(t+1)})|\mathcal{F}^{(t+1)}\right]\right]_{0}\right]\leq\frac{\sigma}{t^{3/2}}.\label{eq:step2-1}
\end{equation}
Summing (\ref{eq:step2-1}) over $t$, we obtain
\[
\sum_{t=1}^{\infty}\mathbb{E}\left[\left[\mathbb{E}\left[L^{(t+1)}(\mathbf{x}^{(t+2)})-L^{(t)}(\mathbf{x}^{(t+1)})|\mathcal{F}^{(t+1)}\right]\right]_{0}\right]<\infty.
\]
Then it follows from the quasi-martingale convergence theorem (cf.
\cite[Th. 6]{Mairal2010}) that $\left\{ L^{(t)}(\mathbf{x}^{(t+1)})\right\} $
converges almost surely.

\textbf{Part 3)} Combining (\ref{eq:basic-relation}) and (\ref{eq:descent-property-1}),
we have
\begin{align}
L^{(t)}(\mathbf{x}^{(t+1)})-L^{(t-1)}(\mathbf{x}^{(t)}) & \leq\nonumber \\
-\eta\bigl\Vert\hat{\mathbf{x}}^{(t)}-\mathbf{x}^{(t)}\bigr\Vert_{2}^{2}+ & L^{(t)}(\mathbf{x}^{(t)})-L^{(t-1)}(\mathbf{x}^{(t)}).\label{eq:step3-2}
\end{align}
Besides, it follows from the convergence of $\left\{ L^{(t)}(\mathbf{x}^{(t+1)})\right\} _{t}$
\[
\lim_{t\rightarrow\infty}L^{(t)}(\mathbf{x}^{(t+1)})-L^{(t-1)}(\mathbf{x}^{(t)})=0,
\]
and the strong law of large numbers that
\[
\lim_{t\rightarrow\infty}L^{(t)}(\mathbf{x}^{(t)})-L^{(t-1)}(\mathbf{x}^{(t)})=0.
\]
Taking the limit inferior of both sides of (\ref{eq:step3-2}), we
have
\[
\begin{split}0 & =\underset{t\rightarrow\infty}{\lim\inf}\left\{ L^{(t)}(\mathbf{x}^{(t+1)})-L^{(t-1)}(\mathbf{x}^{(t)})\right\} \\
 & \leq\underset{t\rightarrow\infty}{\lim\inf}\left\{ -\eta\bigl\Vert\hat{\mathbf{x}}^{(t)}-\mathbf{x}^{(t)}\bigr\Vert_{2}^{2}+L^{(t)}(\mathbf{x}^{(t)})-L^{(t-1)}(\mathbf{x}^{(t)})\right\} \\
 & \leq\underset{t\rightarrow\infty}{\lim\inf}\left\{ -\eta\bigl\Vert\hat{\mathbf{x}}^{(t)}-\mathbf{x}^{(t)}\bigr\Vert_{2}^{2}\right\} \\
 & \qquad\quad+\underset{t\rightarrow\infty}{\lim\sup}\left\{ L^{(t)}(\mathbf{x}^{(t)})-L^{(t-1)}(\mathbf{x}^{(t)})\right\} \\
 & =-\eta\cdot\underset{t\rightarrow\infty}{\lim\sup}\bigl\Vert\hat{\mathbf{x}}^{(t)}-\mathbf{x}^{(t)}\bigr\Vert_{2}^{2}\leq0,
\end{split}
\]
so we can infer that $\lim\sup_{t\rightarrow\infty}\bigl\Vert\hat{\mathbf{x}}^{(t)}-\mathbf{x}^{(t)}\bigr\Vert_{2}=0$.
Since $0\leq\lim\inf_{t\rightarrow\infty}\bigl\Vert\hat{\mathbf{x}}^{(t)}-\mathbf{x}^{(t)}\bigr\Vert_{2}\leq\lim\sup_{t\rightarrow\infty}\bigl\Vert\hat{\mathbf{x}}^{(t)}-\mathbf{x}^{(t)}\bigr\Vert_{2}=0$,
we can infer that $\lim\inf_{t\rightarrow\infty}\bigl\Vert\hat{\mathbf{x}}^{(t)}-\mathbf{x}^{(t)}\bigr\Vert=0$
and thus $\lim_{t\rightarrow\infty}\bigl\Vert\hat{\mathbf{x}}^{(t)}-\mathbf{x}^{(t)}\bigr\Vert=0$.

Consider any limit point of the sequence $\left\{ \mathbf{x}^{(t)}\right\} _{t}$,
denoted as $\mathbf{x}^{(\infty)}$. Since $\hat{\mathbf{x}}$ is
a continuous function of $\mathbf{x}$ in view of (\ref{eq:x-hat-scalar})
and $\lim_{t\rightarrow\infty}\left\Vert \hat{\mathbf{x}}^{(t)}-\mathbf{x}^{(t)}\right\Vert _{2}=0$,
it must be $\lim_{t\rightarrow\infty}\hat{\mathbf{x}}^{(t)}=\hat{\mathbf{x}}^{(\infty)}=\mathbf{x}^{(\infty)}$,
and the minimum principle in (\ref{eq:minimum-principle}) can be
simplified as
\[
(x_{k}-x_{k}^{(\infty)})(\nabla_{k}f^{(\infty)}(\mathbf{x}^{(\infty)})+\xi_{k}^{(\infty)})\geq0,\;\forall x_{k},
\]
whose summation over $k=1,\ldots,K$ leads to
\[
(\mathbf{x}-\mathbf{x}^{(\infty)})^{T}(\nabla f^{(\infty)}(\mathbf{x}^{(\infty)})+\boldsymbol{\xi}^{(\infty)})\geq0,\quad\forall\mathbf{x}.
\]
Therefore $\mathbf{x}^{(\infty)}$ minimizes $L^{(\infty)}(\mathbf{x})$
and $\mathbf{x}^{(\infty)}=\mathbf{x}^{\star}$ almost surely by Lemma
\ref{lem:connection}. Since $\mathbf{x}^{\star}$ is unique in view
of Assumptions (A1)-(A3), the whole sequence $\{\mathbf{x}^{(t)}\}$
has a unique limit point and it thus converges to $\mathbf{x}^{\star}$.
The proof is thus completed.
\end{IEEEproof}

\section*{Acknowledgment}

The authors would like to thank Prof. Gesualdo Scutari for the helpful
discussion.

\end{document}